\documentclass{amsart}

\usepackage{amssymb}
\usepackage[all]{xy}
\usepackage{hyperref}

\usepackage{enumitem}   %from SK, fixes \enumerate

\setlist[enumerate]{itemsep=.2em,topsep=.2em,leftmargin=1.25em,itemindent=2.0em}

\newtheorem{thm}{Theorem}%[section]

\newtheorem{lem}[thm]{Lemma}
\newtheorem{cor}[thm]{Corollary}

\newtheorem{prop}[thm]{Proposition}

\newtheorem{complem}[thm]{Complement}%!!!!!!!!!!!!!!!!!!!!!!
\theoremstyle{definition}
\newtheorem{defn}[thm]{Definition}

\newtheorem{say}[thm]{}
\newtheorem{exmp}[thm]{Example}

\newtheorem{rem}[thm]{Remark}          

\newtheorem*{ack}{Acknowledgments}      % \renewcommand{\theack}{} 
\newtheorem{notation}[thm]{Notation}   
  
\newtheorem{defn-thm}[thm]{Definition--Theorem}  %!!!!!!!!!!!!!!!!!!!!!!!!
\newtheorem{defn-lem}[thm]{Definition--Lemma}  %!!!!!!!!!!!!!!!!!!!!!!!!
\newtheorem{complem-p}[thm]{Complement}

\theoremstyle{remark}

\renewcommand{\o}[0]{{\mathcal O}} 
\newcommand{\z}[0]{{\mathbb Z}}
\newcommand{\n}[0]{{\mathbb N}}
\renewcommand{\r}[0]{{\mathbb R}} 

\renewcommand{\a}[0]{{\mathbb A}}

\newcommand{\p}[0]{{\mathbb P}}

\newcommand{\q}[0]{{\mathbb Q}}
\newcommand{\map}[0]{\dasharrow}
\newcommand{\qtq}[1]{\quad\mbox{#1}\quad}
\newcommand{\spec}[0]{\operatorname{Spec}}

\newcommand{\gal}[0]{\operatorname{Gal}}

\newcommand{\rank}[0]{\operatorname{rank}}

\newcommand{\supp}[0]{\operatorname{Supp}}    
\newcommand{\red}[0]{\operatorname{red}}

\newcommand{\Hom}[0]{\operatorname{Hom}}

\newcommand{\chow}[0]{\operatorname{Chow}}
\newcommand{\chr}[0]{\operatorname{char}}

\newcommand{\hilb}[0]{\operatorname{Hilb}}

\newcommand{\mor}[0]{\operatorname{Mor}} 
 
\newcommand{\tsum}[0]{\textstyle{\sum}}

\newcommand{\shom}[0]{\operatorname{\mathcal{H}\!\it{om}}}

\newcommand{\ind}[0]{\operatorname{ind}}

\def\into{\DOTSB\lhook\joinrel\to}

\def\loccoh#1.#2.#3.#4.{H^{#1}_{#2}(#3,#4)}

\DeclareMathAlphabet{\mathchanc}{OT1}{pzc}%
                                {m}{it}

\newcommand{\PGL}{\mathrm{PGL}}

\usepackage[all]{xy}\xyoption{dvips}

\newcommand{\hci}[0]{$H^{\rm ci}$\ }

\newcommand{\bideg}[0]{\operatorname{bideg}}

\newcommand{\ecomb}[0]{\operatorname{EComb}}
\newcommand{\env}[0]{\operatorname{Env}}

\newcommand\cupn{\mathbin{\ooalign{$\cup$\cr%
   \hfil\raise0.42ex\hbox{$\scriptscriptstyle\times$}\hfil\cr}}}

\newcommand\cupf{\mathbin{\ooalign{$\cup$\cr%
   \hfil\raise0.42ex\hbox{$\scriptscriptstyle f$}\hfil\cr}}}

\begin{document}
\bibliographystyle{amsalpha}

%% \hfill\today

\title[Stable maps of curves]{Stable maps of curves and \\ algebraic equivalence of 1-cycles}
        \author{J\'anos Koll\'ar and Zhiyu Tian}

        \begin{abstract}  We show that algebraic equivalence of images of stable maps of curves lifts to deformation equivalence of the stable maps.
          The main applications concern $A_1(X)$,  the group of 1-cycles modulo algebraic equivalence, for smooth, separably rationally connected varieties.
          If $K/k$ is an algebraic extension, then
          the kernel of $A_1(X_k)\to A_1(X_K)$ is at most $\z/2\z$.
          If $k$ is finite, then the image equals  the subgroup of Galois invariant cycles.

          %% This paper replaces Sections~2--3 of 2211.15915.v.1 and  Sections~2--3 of 2211.15911.v.1. The other Sections are retained in the revised versions of these papers. 
          
        \end{abstract}

\maketitle

For a scheme $X$, let $\mor(\mbox{Curves}, X)$ denote the stack of morphisms from proper, nodal curves to $X$. If two morphisms $\pi_i:C_i\to X$ are
deformation equivalent (see Definition~\ref{def.curves.defn}), then
the image 1-cycles  $(\pi_i)_*[C_i]$ are algebraically equivalent. The converse, however,
rarely holds.

Nonetheless, we show that if $X$ is smooth, then algebraic equivalence of the image cycles lifts to deformation equivalence of morphisms.
We use $A\cupn B$ to indicate that all intersection points $A\cap B$ are nodes; see Notation~\ref{amal.cupn.not}.

  \begin{thm} \label{main.same.tail.thm}
    Let $X$ be a smooth, projective  variety over an algebraically closed  field $K$. Let 
    $\pi_i:C_i\to X$  (for $i\in I$) be finitely many morphisms of nodal curves   to $X$ such that the 
    $(\pi_i)_*[C_i]$ are  algebraically equivalent to each other.
    Then there is a  nodal  deformation equivalence
    (as in Definition~\ref{def.curves.defn})
  $$
\begin{array}{ccccc}
C_i\cupn  R_i & \subset  & S  &  \stackrel{\pi}{\longrightarrow} & B\times X\\
  \downarrow &  &  \downarrow &  &  \downarrow \\
  b_i  & \in  & B   & = & B,
\end{array}
\eqno{(\ref{main.same.tail.thm}.1)}
$$
such that
\begin{enumerate}\setcounter{enumi}{1}
  \item $(\pi|_{C_i}:C_i\to X) \cong (\pi_i:C_i\to X)$  for $i\in I$, and the
  \item $(\pi|_{R_i}: R_i\to X)$ are isomorphic to each other  for $i\in I$.
\end{enumerate}
  \end{thm}

 That is, the
 algebraic equivalence is visible at the level of maps, not just at the level of cycles.

 For arithmetic applications, the following version is quite useful.

 \begin{thm} \label{main.same.tail.thm.k}
  Let $k$ be a perfect field with algebraic closure $K$, and $k\subset L\subset K$.
    Let $X$ be a smooth, projective, $k$-variety
  and  $\pi_i:C_i\to X_L$ (for $i\in I$) morphisms of nodal curves defined over $L$,  such that the 
    $(\pi_i)_*[C_i]$ are  algebraically equivalent to each  other over $K$.
Assume also that $\chi(C_i,\o_{C_i})$ is independent of $i\in I$.

Then one  can choose  (\ref{main.same.tail.thm}.1) such that, in addition to
(\ref{main.same.tail.thm}.2--3), 
\begin{enumerate}
\item the $\pi|_{C_i\cupn R_i}:C_i\cupn R_i\to X_L$ are  defined over $L$  for $i\in I$, and the
\item $(\pi|_{R_i}: R_i\to X_L)$ are defined over $k$.
   \end{enumerate}
\end{thm}
  
 \begin{rem}\label{main.same.tail.rem} There are some unusual features of
Theorem~\ref{main.same.tail.thm.k}.

    First, we  assume only that   the 
    $(\pi_i)_*[C_i]$ are  algebraically equivalent over $K$. Correspondingly, the diagram
    (\ref{main.same.tail.thm}.1) itself is defined only over $K$.
    It is only the fibers over the points $b_i$ that are defined over  $L$.

    Second, the $R_i$ are defined over the smallest possible field, $k$.

    Third, the assumption  on $\chi(C_i,\o_{C_i})$ is necessary, see Example~\ref{quad.deg2.exmp}. 
 \end{rem}

 From Theorem~\ref{main.same.tail.thm.k}, we can not conclude that the $C_i\cupn R_i$ are algebraically equivalent over $k$ or $L$, since the
 families $S\to B$ are only defined over $K$. However, we will be able to prove
 algebraic equivalence for certain varieties, which we define next.

\begin{defn} \label{src.cd.1.defn}  Let $X$ be a  smooth, proper variety  over an  algebraically closed field. Given any $g:\p^1\to X$, the pull-back 
  $g^*T_X$ decomposes as $ \sum_i \o_{\p^1}(a_i)$.  It is called {\it very free} if $a_i>0$ for every $i$, and  {\it almost very free} if $a_i\geq 0$ for every $i$, with strict inequality for all but one $a_i$.

  $X$ is 
 {\it separably rationally connected,} or {\it src,} if there is a very free $g$; see \cite[Sec.IV.3]{rc-book} or
 \cite[Sec.7]{ar-ko} for many equivalent definitions.
 
 We say that $X$ is {\it separably rationally connected in codimension 1,} or {\it src in codimension 1,} if there is an  almost very free $g:\p^1\to X$.  Thus src implies src in codimension 1. The main new examples are  varieties that  admit  a morphism to a curve  $X\to B$ whose general fibers are smooth and   separably rationally connected.
 These 2 classes give  all examples in characteristic 0.
 
  There are a few more  in positive characteristic.
  General fibers of  $X\to B$ could be  singular,  but contain very free  rational curves in their smooth locus. There are also rationally connected (but not src) examples, studied in   \cite{MR2629599}.

 Being src (or src in codimension 1) are both  open conditions in smooth families of projective varieties.
  \end{defn}

For varieties that are src in codimension 1, we  turn
the deformation  equivalences in  Theorem~\ref{main.same.tail.thm.k}
into algebraic equivalences over $L$.  This leads to a series of arithmetic consequences for 1-cycles.
Further applications to  the Brauer-Manin obstruction  for zero cycles on  geometrically rational surfaces  over  global function fields and the coniveau filtration  are discussed in
\cite{https://doi.org/10.48550/arxiv.2211.15915, https://doi.org/10.48550/arxiv.2211.15911}.
\medskip

\begin{say}[Algebraic equivalence of 1-cycles] \label{al.eq.1c.say}
  Let $X$ be a proper scheme over a field $k$. We use  $A_d(X)$ for the group of $d$-dimensional cycles modulo {\it algebraic equivalence,} the latter denoted by $Z_1\sim_a Z_2$. (For an introduction, see
  \cite[Sec.10.3]{MR1644323}, where this group is denoted by $B_d(X)$.)

  Let $L/k$ be an algebraic  field extension.  We have natural maps
  $A_d(X_k)\to A_d(X_L)$. The kernel of this is $\deg(L/k)$-torsion.

  If $L/k$ is Galois, then the image is contained in
   the Galois-invariant subgroup $A_d(X_L)^{\gal(L/k)}$. Thus our interest is the
  kernel and cokernel of the map
  $$
  A_1(X_k)\to A_1(X_L)^{\gal(L/k)}.
  \eqno{(\ref{al.eq.1c.say}.1)}
    $$
\end{say}

The following 3 theorems are proved in Section~\ref{arith.sec}.
The first one  describes  the kernel of (\ref{al.eq.1c.say}.1).

\begin{thm}\label{main.thm.1}
  Let $X_k$ be a smooth, projective   variety over  a perfect field $k$ with algebraic closure $K$. Assume that $X_K$ is  src in codimension 1.
    Then the kernel of the  natural map
  $A_1(X_k)\to A_1(X_K)$ is either trivial or $\z/2\z$.
 More precisely,
  \begin{enumerate} %%\setcounter{enumi}{2}
    \item the kernel is
      trivial if  $X_k$ contains an odd degree 0-cycle, and
      \item if  $Z=\sum d_iC_i$  and $Z_K\sim_a 0$, then $Z\sim_a 0$ iff 
  the index of $X$ (Definition~\ref{ind.ech.defn})
   divides $\chi(Z):=\sum_i d_i\chi(C_i, \o_{C_i})$.
\end{enumerate}
\end{thm}

 For finite fields (\ref{al.eq.1c.say}.1) is an isomorphism.

\begin{thm}\label{main.k.src.thm.cor.3} Let $k$ be  a  perfect field with algebraic closure $K$, and
  $X$  a smooth, projective $k$-variety that is  src in codimension 1.
  Assume that every geometrically irreducible $k$-variety has a 0-cycle of degree 1. (For example, $k$ is finite or pseudo algebraically closed.)
  Then
  $$
  A_1(X_k)\to A_1(X_K)^{\gal(K/k)}\qtq{is an isomorphism.}
  \eqno{(\ref{main.k.src.thm.cor.3}.1)}
    $$
\end{thm}

Let $k$ be  a field, $C$ a smooth, projective, geometrically irreducible $k$-curve. Then $X:=C\times \p^1$ is  src in codimension 1. For $c\in C(K)$, the line
$\{c\}\times \p^1$ gives a class in $A_1(X_K)^{\gal(K/k)}$.
If (\ref{main.k.src.thm.cor.3}.1) holds then $C$ has a 0-cycle of degree 1.
So the assumption on $k$ in Theorem~\ref{main.k.src.thm.cor.3} is necessary in all cases. 
\medskip

The third application concerns specialization of 1-cycles.

\begin{thm}\label{main.k.src.thm.cor.2} Let $R$ be a Henselian DVR  with perfect residue field. Let $p:X\to \spec R$ be a smooth, projective morphism with
  closed fiber $X_0$ and generic fiber $X_g$.
  Assume that $X_0$  is src in codimension 1.

  Then the specialization map
  $A_1(X_g)\to A_1(X_0)$ is an isomorphism.
\end{thm}

Next we give some examples concerning Theorems~\ref{main.same.tail.thm.k} and \ref{main.thm.1}.
From now on we write $\chi(C):=\chi(C, \o_C)$ for proper, 1-dimensional schemes.

\begin{exmp} \label{quad.deg2.exmp}
Take $X=\p^3$, choose $m>0$ and let $C_1\subset \p^3$ be a degree $2m$ curve of genus 1 over $\q$ such that  every 0-cycle has degree divisible by $2m$.
    Let $C_2$ be the union of $m$ disjoint, conjugate conics, such that again 
    every 0-cycle has degree divisible by $2m$. Note that
    $\chi(C_1)=0$ and $\chi(C_2)=m$. Let $R$ be any curve and attach it to $C_i$ at $n_i$ nodes, where, necessarily,   $2m\mid n_i$. Since
    $\chi(C_i\cupn R)=\chi(C_i)+\chi(R)+n_i$, we obtain that
    $$
    \chi(C_1\cupn R)-\chi(C_2\cupn R)\equiv m\mod 2m.
    $$
    The Euler characteristic is deformation invariant, so
    the $C_i\cupn R$ are not deformation equivalent.
\end{exmp}
 
\begin{exmp} \label{first.exmaples.exmp}
  In connection with Theorem~\ref{main.thm.1},
we give 3 examples computing the kernel of $A_1(X_k)\to A_1(X_K)$.

 (\ref{first.exmaples.exmp}.1)  Let $Q_k\subset \p^4$ be a quadric without $k$-points.
  Let $C_1\subset Q_k$ be a smooth conic and $C_2$ the union of a conjugate pair of lines.
  For example, for any $k\subset \r$ we can take
  $$
  \begin{array}{l}
  Q=(\tsum_i x_i^2=0)\subset \p^4,\ C_1=(x_3=x_4=0), \qtq{and}\\
  C_2=(x_0=x_1-ix_2=x_3-ix_4=0)\cup (x_0=x_1+ix_2=x_3+ix_4=0).
  \end{array}
$$
Then $C_1-C_2$ is  algebraically equivalent to 0 over $K$, but 
not algebraically equivalent to 0 over $k$ by
 (\ref{chi.alg.eq.lem})
since $\chi(C_1)=1$ and $\chi(C_2)=2$.

Thus the kernel of  $A_1(Q_k)\to A_1(Q_K)$ is $\z/2\z$.

(\ref{first.exmaples.exmp}.2)
 Let $\pi:X\to B$ be a morphism of smooth projective varieties.
 Assume that $\dim B=1$ and the generic fiber has odd index.
 Write the index of $B$ as $2^a(\mbox{odd})$. Then  the index of $X$ is also $2^a(\mbox{odd})$.  Let $Z_1, Z_2$ be 1-cycles on $X$ that are 
 algebraically equivalent  over $K$. Then the $Z_i\to B$ have the same degree, so 
 $\chi(Z_1)\equiv \chi(Z_2)\mod 2^a$ by (\ref{chi.alg.eq.lem}.1).
 Using (\ref{ind.ech.defn}) we get that the index of $X$ divides
 $\chi(Z_1)- \chi(Z_2)$.

 So, if the generic fiber is src, then
 $A_1(X_k)\to A_1(X_K)$ is an injection.

 (\ref{first.exmaples.exmp}.3)
 Let $\pi:X\to S$ be a morphism of smooth projective varieties.
 Assume that $\dim S=2$ and the generic fiber  has odd index.
 As before, if $\ind_S=2^a(\mbox{odd})$, then   $\ind_X=2^a(\mbox{odd})$.  Let $Z_1, Z_2$ be 1-cycles on $X$ that are 
 algebraically equivalent  over $K$. Then $\pi_*Z_1$ and
 $\pi_*Z_2$ are also algebraically equivalent  over $K$. By the adjunction formula, they have the same Euler characteristic. This implies that
 $\chi(Z_1)\equiv \chi(Z_2)\mod 2^a$ by (\ref{chi.alg.eq.lem}.1).
 Using (\ref{ind.ech.defn}) we get that the index of $X$ divides
 $\chi(Z_1)- \chi(Z_2)$.

 So, if $S$ is geometrically rational or ruled, and the generic fiber is src,   then
 $A_1(X_k)\to A_1(X_K)$ is an injection.
\end{exmp}

\begin{say}[Steps of the proof]\label{steps.pf.say}

   Assume for simplicity that $C_1, C_2\subset X$ are algebraically equivalent smooth curves of the same genus. By definition, their algebraic equivalence is certified (\ref{certifies.defn}) by a
deformation
   $$
  \begin{array}{ccccc}
   S_i & \subset  & S  &  \stackrel{\pi}{\longrightarrow} & B\times X\\
  \downarrow &  &  \downarrow &  &  \downarrow \\
  b_i  & \in  & B   & = & B,
  \end{array}
  \eqno{(\ref{steps.pf.say}.1)}
  $$
  such that  $\pi_*[S_1]=[C_1]+Z$ and $\pi_*[S_2]=[C_2]+Z$ for some
  1-cycle $Z$. After applying semistable reduction to $S\to B$, we may assume that all fibers of $S\to B$ are nodal. 

  Note that in general the $S_i$ need not have irreducible components $C'_i\subset S_i$ that are isomorphic to the $C_i$. Even if they do, and we write $S_i=C'_i\cup R_i$, the curves $R_i$ and the maps  $\pi_i:=\pi|_{R_i}:R_i\to \supp Z$ can be
  very different from each other.

  The moduli space of morphisms of curves to $\p^1$ was studied by 
   Hurwitz. Similarly, we  view the $R_i$ as `covers' of the same $\supp Z$, though the maps $R_i\to\supp Z$ need not even be finite. 
  In Section~\ref{map.cur.sec} we give a combinatorial description of going between various such `covers' by gluing  additional  components and deformation equivalences.

  Next we want to realize the above combinatorial steps as deformation equivalences. The general problem is essentially the following:

  \medskip
 {\it Question \ref{steps.pf.say}.2.}  Write a reducible curve as $C_1\cup C_2$ and construct deformations of
 each $C_i$. Can they be glued to a deformation of  $C_1\cup C_2$?
 \medskip

 If the $C_i$ are disjoint then yes, but  we run into non-flatness at the intersection points  $C_1\cap C_2$. However, flatness can be restored by adding an extra component. The general framework is discussed in
 Section~\ref{glue.sec} and applied to our current setting in 
 Section~\ref{modify.sec}.
The main technical  result is
 (\ref{tail.same.thm}), 
  which replaces (\ref{steps.pf.say}.1)
  with another certificate
   $$
  \begin{array}{ccccc}
   C_i\cupn R_i & \subset  & S_R  &  \stackrel{\pi_R}{\longrightarrow} & B\times X\\
  \downarrow &  &  \downarrow &  &  \downarrow \\
  b_i  & \in  & B   & = & B,
  \end{array}
  \eqno{(\ref{steps.pf.say}.3)}
  $$
  such that now  $(\pi_R:R_1\to X)\cong (\pi_R:R_2\to X)$. (\ref{tail.same.thm}) also gives very detailed information on how $S_R$ is obtained from $S$.

  This completes the  geometric part of the argument, proving
  Theorem~\ref{main.same.tail.thm}.

  If the $C_i$ are defined over some field $k$, then we first construct
  (\ref{steps.pf.say}.3) over the algebraic closure  $K$, and then in
  (\ref{main.same.tail.thm.p.4}) construct  another certificate,
  where 
  $$
  C_i+R_i\qtq{is replaced by} C_i+(\mbox{all Galois conjugates of } R_i).
  $$
  This gives us  Theorem~\ref{main.same.tail.thm.k} and completes the  part of the proof that holds for all smooth, projective varieties.

  The base curve $B$ of the deformation in (\ref{steps.pf.say}.3) is
  geometrically connected, but usually not definable over $k$, so it says nothing about algebraic equivalence over $k$.

  Now we face the following
  
  \medskip
      {\it Question \ref{steps.pf.say}.4.}  Let $Y$ be a $k$-variety, and $Z_i\subset Y$ curves such that $[Z_1]$ and $[Z_2]$  are in the same  geometrically connected component of $\hilb_1(Y)$. Are $[Z_1]$ and $[Z_2]$ algebraically equivalent over $k$?
 \medskip

 The answer is negative in general, and
  this is where we  need that $X$ be src in codimension 1.
    We show in (\ref{main.src.comb.thm.c}) that there is a third curve $T$ such that
  $[Z_1\cup T]$ and $[Z_2\cup T]$ are smooth points on the same
   irreducible component of  $\hilb_1(Y)$ that is also geometrically irreducible. Thus
  $[Z_1\cup T ]$ and $[Z_2\cup T]$ are algebraically equivalent over $k$, and
   so are $[Z_1]$ and $[Z_2]$.

   This gives us 
   Theorem~\ref{main.thm.1}, while   Theorems~\ref{main.k.src.thm.cor.3}--\ref{main.k.src.thm.cor.2} require a more careful understanding of the above procedure.
\end{say}

\begin{ack}  We thank  Jean-Louis~Colliot-Th\'el\`ene, Rahul~Pandharipande and Olivier~Wittenberg for many helpful and constructive comments. 
Partial  financial support  to JK  was provided  by  the NSF under grant number
DMS-1901855, and to ZT  by NSFC grants No.11890660 and No.11890662.
 \end{ack}

 \section{Notation and definitions}

\begin{defn} \label{def.curves.defn}
  Let $X$ be a proper scheme over an algebraically closed  field $K$.
A {\it deformation equivalence} (of curves mapping to $X$) is a diagram
   $$
  \begin{array}{ccccc}
   S_i & \subset  & S  &  \stackrel{\pi}{\longrightarrow} & B\times X\\
  \downarrow &  &  \downarrow &  &  \downarrow \\
  b_i  & \in  & B   & = & B,
  \end{array}
  \eqno{(\ref{def.curves.defn}.1)}
  $$
  such that
  \begin{enumerate}\setcounter{enumi}{1}
\item $B$ is a  connected curve with smooth marked points $b_i$, 
\item  $S\to B$ is flat, proper,  of pure relative dimension 1, and
   \item   $S_i:=S_{b_i}$ is the fiber over $b_i$.
  \end{enumerate}
  The $b_i$ and the intersection points of different irreducible components of $B$ are the {\it pivot points.}
  
  Such a deformation is called  {\it nodal} (or 2-nodal) if all fibers are nodal curves.

   We call a deformation {\it 3-nodal} if all but finitely many fibers are nodal curves, and the pivot fibers have  only  nodes  (like 2 coordinate axes  $(xy=0)\subset \a^2$)   and triple points (like 3 coordinate axes 
   $(xy=yz=zx=0)\subset \a^3$).

   Our constructions naturally produce 3-nodal deformations, we then make them 2-nodal at the end using (\ref{jump.over.say}). 

  {\it Note on terminology.} 3-nodal is not standard. These are sometimes called ordinary triple points, but that name usually refers to planar triple points whose tangent cone is 3 distinct lines.
\end{defn}

\begin{say}[Base change]
  Let  $(B', b'_i)$ be another    connected curve with smooth marked points,  and $r: B'\to B$ a morphism such that $r(b'_i)=b_i$. By pull-back we get a
  new deformation equivalence
  $$
  B'\leftarrow S':=B'\times_BS \to B'\times X.
  $$
  If $S\to B$ is nodal, then so is $S'\to B'$. If $S'\to B'$ is 3-nodal,
  and the fibers of $S\to B$ over the images of the pivot points of $B'$ are 3-nodal, then $S'\to B'$ is also 3-nodal. (For us $B'\to B$ will usually come from a general choice, so the extra condition will not be a problem.)

  During the proofs we have to make many base changes, and usually we use $B$ to denote all base curves. 
  \end{say}

\begin{say}[Improving the singularities of $S\to B$]\label{imp.sing.s.b.say}
  Let $p:S\to B$ be a flat, projective  morphism of a surface to an irreducible curves. The semistable reduction theorem \cite{del-mum} says that there is a base change
  $B^{\rm ss}\to B$, and a resolution of singularities
  $S^{\rm ss}\to B^{\rm ss}\times_BS$   such that the resulting
  $S^{\rm ss}\to B^{\rm ss}$ is smooth over generic points, and the singular fibers are all nodal.

  If $B$ is reducible, we could apply this to each irreducible component
  $S_i\to B_i$. However, the fibers of
  $S^{\rm ss}_i\to B^{\rm ss}_i$ and of $S^{\rm ss}_j\to B^{\rm ss}_j$ over $B_i\cap B_j$ are usually not isomorphic. So the families $S^{\rm ss}_i\to B^{\rm ss}_i$ can not be glued together into a flat family over a connected curve.

  Nonetheless, if we have a  deformation equivalence
  as in (\ref{def.curves.defn}.1) where all but finitely many fibers are nodal,
  then, after a suitable base change  $B'\to B$, we can get a new family
  $S'\to B'$ where all fibers over smooth points of $B'$ are nodal and
  the fibers over the pivot points are unchanged.

  For this reason, we can and will ignore possible badly singular fibers over 
  smooth points of $B$.

  \end{say}

\begin{defn} \label{certifies.defn}
  Let $X$ be a proper scheme over a field $k$ and $Z$ a 1-cycle on $X$. We say that  the deformation (\ref{def.curves.defn}.1)
    {\it certifies} that $Z\sim_a 0$ iff
    $\pi_*[S_1]-\pi_*[S_2]=Z$, as 1-cycles on $X$.

    By definition, if $Z\sim_a 0$, there is always a certificate
    $B\leftarrow  S \to B\times X$ where $B$ is a smooth,  irreducible curve.
    We can apply stable reduction to $B\leftarrow  S$ and get
    certificate
    $B\leftarrow  S^{\rm nodal} \to B\times X$ where the fibers of
    $S^{\rm nodal}\to B$ are nodal curves.

    We can also assume that the general fibers are smooth, but this usually makes the $S_i$ disconnected. For us, connectedness is more  important.
     \end{defn}

\begin{notation} \label{amal.cupn.not}
  Let $X_i$ be schemes, $Z_i\subset X_i$ subschemes and
  $\sigma:Z_1\cong Z_2$ an isomorphism. We use
  $X_1\amalg_\sigma X_2$  to denote the scheme obtained by identifying
  $Z_1$ and $Z_2$.

  Let $C_i$ be curves, $Z_i\subset C_i$  sets of smooth points  and
  $\sigma:Z_1\cong Z_2$ an isomorphism.  Instead of $C_1\amalg_\sigma C_2$ we use 
  $C_1\cupn_{\sigma} C_2$ (where  the $\times$ reminds us that we get a node). We use  $C_1\cupn C_2$
  if the $Z_i$ and $\sigma$ are clear or not important.
    \end{notation}

\begin{say}[From 3-nodes to  nodes]\label{jump.over.say}
  Let $C_3\subset \p^3$ be  3 general lines meeting at point $p$.
  Using (\ref{hilb.defn}.1--2) we compute that the Hilbert scheme of $\p^3$ is smooth of dimension 12 at $[C_3]$, and includes the (3-dimensional) universal deformation of the singularity at the intersection point.

  Blow up $p$, let $E\cong \p^2$ be the exceptional divisor.  Let $\bar C_3$ be a curve consisting of  the birational transform $C'_3$ of $C_3$, and a conic  $C_0\subset E$ passing through the 3 points  $E\cap C'_3$. 
  Although the conic is not unique, 
  $\pi:\bar C_3\to C_3$ is unique as an abstract curve.

  The space of morphisms  $\mor(\mbox{Curves}, \p^3)$ is also
  smooth of dimension 12 at $[\pi:\bar C_3\to \p^3]$; cf.\ \cite[15]{ar-ko}.
  Since $R^1\pi_*\o_{\bar C_3}=0$, $\pi$ extends to a morphism
  $$
  \Pi: \mor(\mbox{Nodal curves}, \p^3)^\circ
  \to \hilb_1(\p^3)^\circ
  \eqno{(\ref{jump.over.say}.1)}
  $$
  in suitable open neighborhoods of $[\pi:\bar C_3\to \p^3]$ and  $[C_3]$.
  Thus $\Pi$ is a local isomorphism. That is, every flat deformation of the 3-node
  $C_3$ can be lifted to a flat deformation of the morphism
  $\pi:\bar C_3\to C_3$. Since $\bar C_3$ has only nodes,
  this replaces 3-nodal families with nodal ones.
In general, for any $X$ we get 
$$
  \Pi_X: \mor(\mbox{Nodal curves}, X)
  \map \hilb_1(X),
  \eqno{(\ref{jump.over.say}.2)}
  $$
  which restricts to an  isomorphism over the 3-nodal locus of $\hilb_1(X)$.
  \smallskip
  
  (In fact, nodes and 3-nodes are the only singularities over which
   $\Pi_X$ is a local isomorphism.)
\smallskip

  The following explicit example may be useful.
It shows how to    `jump over' the  nodes
as we slide  an auxiliary curve $L$ from one irreducible component of a nodal curve $C$ to another.

  Consider a family of nodes parametrized by a node:
 $$
 \a^4_{xyuv}\supset S:=(xy=uv=0) \stackrel{\pi}{\longrightarrow}
 (xy=0)=:B \subset \a^2_{xy}.
 $$
 Let $D\subset S$ be the image of the
  section  $\sigma(x,y)=(x,y,x,y)$. Thus, along the $x$-axis we move in the
 $u$-axis, along the $y$-axis we are in the
 $v$-axis.  At the origin  $(x=y=0)$ we hit the node $(u=v=0)$.

Let $D'\subset S$ be the image of the
section $\sigma'(x,y)=(x,y,y,x)$. Its equation is $(v-x=u-y=0)$.

Note that $D+D'$ is  Cartier divisor,
with defining equation $x+y=u+v$.

Blow  up  $D'\subset S$ to get $\bar S\to S$. The restriction of $D'$ to the
$(x=v=0)$ and  $(y=u=0)$ planes is a Cartier divisor, so on these the blow-up is the identity. The restriction of $D'$ to the
$(x=u=0)$ and  $(y=v=0)$ planes is the origin, so we get their usual blow-up.
The 2 exceptional curves are identified in $\bar S$ to get $\bar E\subset \bar S$.

The new central fiber consists of $\bar E$ and the birational transforms of
$(0,0,u,0)$ and $(0,0,0,v)$, meeting $\bar E$ at distinct points $p_u, p_v$.
The birational transform of $D$ is still a section $\bar D$, meeting
$\bar E$ at a point different from the $p_u, p_v$. Thus in the new family
$\bar S\to B$, the  section $\bar D$ is a  Cartier divisor that does not pass through any of the nodes of the fibers.
\end{say}

  \section{Maps to curves}\label{map.cur.sec}

Hurwitz proved that for fixed $(g, d)$, all degree $d$ morphisms from  genus $g$ smooth curves to $\p^1$ form an irreducible family.
We discuss what happens when $\p^1$ is replaced by an arbitrary (possibly  reducible) curve  $C$ and we study maps from nodal curves to $C$.

Irreducibility fails already in simple cases. For example, let
$C:=(xyz=0)\subset \p^2$ be a cycle of 3 lines. Separating one of its nodes we get 3 morphisms  $\pi_i:A\to C$ where $A$ is a chain of 3 rational curves, so $\chi(A)=1$.
These 3 maps are distinct, rigid and all their numerical invariants agree.

We aim to prove that by adding more irreducible curves to $A$, we get a connected moduli space.

\begin{defn} \label{elem.ch.over.c.def}
  Let $k$ be an algebraically closed field, $C$ a proper curve, $A, B$ proper nodal curves and $\pi_A:A\to C,\pi_B:B\to C$ morphisms.
  The following operations are the  {\it elementary changes} of these data.
  \begin{enumerate}
  \item Take a smooth, projective, irreducible curve  $\tau:E\to C$ and attach it to $A$ and $B$ along  two subsets of the same cardinality, to get
    $\pi'_A:A':=A\cupn E\to C$ and $\pi'_B:B':=B\cupn E\to C$.
    \item Take constant morphisms $\tau_A:\p^1\to C, \tau_B:\p^1\to C$ and 
      attach then to $A$ (resp.\  $B$) along  two subsets of the same cardinality, to get
      $\pi'_A:A':=A\cupn \p^1\to C$ and $ \pi'_B:B':=B\cupn \p^1\to C$.
  \item     Write $A=E\cap F$. Choose $\pi'_F:F'\to C$ to be   deformation equivalent to $\pi|_{F}$ fixing $E\cap F$. Set
      $$
      \pi'_A:(\pi|_{E}, \pi'_F): A':=E\cupn F'\to C
      \qtq{and} \pi'_B=\pi_B: B\to C.
      $$
      \item As in (3) with the roles of $A, B$ reversed.
      \end{enumerate}
\end{defn}

The main result of this section is the following.

  \begin{prop}  \label{to.curves.prop}
    Let $k$ be an algebraically closed field, $C$ a proper curve, $A, B$ proper nodal curves and $\pi^1_A:A^1\to C,\pi^1_B:B^1\to C$  morphisms.
    The following are equivalent.
\begin{enumerate}
\item    $(\pi^1_A)_*[A^1]=(\pi^1_B)_*[B^1]$ and $\chi(A^1)=\chi(B^1)$.
\item There is a sequence of elementary changes  (\ref{elem.ch.over.c.def}.1--4)
  $$
  \bigl(\pi^1_A:A^1\to C,\pi^1_B:B^1\to C\bigr)  \map
  \cdots \map
  \bigl(\pi^r_A:A^r\to C,\pi^r_B:B^r\to C\bigr)
  $$
  such that  $\bigl(\pi^r_A:A^r\to C\bigr) \cong\bigl(\pi^r_B:B^r\to C\bigr)$.
    \end{enumerate}
\end{prop}

Proof. It is clear that elementary changes  preserve both 
 $(\pi^i_A)_*[A^i]-(\pi^i_B)_*[B^i]$ and $\chi(A^i)-\chi(B^i)$, so (2) $\Rightarrow$ (1).

For the converse, we may assume that  $\pi^1_A, \pi^1_B$ are dominant.  
The claim is clear if $A^1$ or $B^1$ are empty.

Let $C_j$ be the  irreducible components of $C$ with normalizations $\bar C_j\to C_j$. Take curves  $\bar A^1_j\to \bar C_j$ and $\bar B^1_j\to \bar C_j$ such that $A^1$  (resp $B^1$) is obtained from
$\amalg_j\bar A^1_j$ (resp.\ $\amalg_j\bar B^1_j$) by gluing some point pairs to nodes.
 These $\bar A^1_j, \bar B^1_j$ are almost unique, except
for irreducible components of $A^1\cup B^1$ that lie over a singular point of $C$.
For these we have finitely many possible choices, we make any one of them.

Our plan is to show that our claim holds each
$\bigl(\pi^1_A:\bar A^1_j\to \bar C_j,\pi^1_B:\bar B^1_j\to \bar C_j\bigr)$.
However, $\chi(\bar A^1_j)$ can be different from $\chi(\bar B^1_j)$.
So first we need a series of elementary changes to ensure that
$\chi(\bar A^i_j)=\chi(\bar B^i_j)$ for some $i\geq 1$.

It will be convenient to add disjoint copies of
$\amalg_j\bar C_j$ 
to both $A^1$ and  $B^1$, in order  to avoid some trivial special cases.

We have to pay special attention to the nodes of $A^1$ and $B^1$ that get separated  when we pass to $\amalg_j\bar A^1_j$ and $\amalg_j\bar B^1_j$. So let  $c\in C$ be a singular point and  $C_j(c)$ the local branches of $C$ through $c$.
Let $N^1_A(c, j_1, j_2)\subset  A^1$ be the
set of all  nodes of $A^1$ over $c$ whose 2 branches end up on
$C_{j_1}(c)$ and  $C_{j_2}(c)$ when we pass to $\amalg_j\bar A^1_j$.
We call these {\it separating nodes.}
We  start with some elementary changes to achieve that 
$\# N^i_A(c, j_1, j_2)=\# N^i_B(c, j_1, j_2)$ for every $c, j_1, j_2$.

So, let us see how to create a new node on $A^1$ by elementary changes. This is done in 2 steps.
\begin{enumerate}
\item First attach a copy of $E_1\cong \bar C_{j_1}$ to both $A^1$ and $B^1$ at a point  as in (\ref{elem.ch.over.c.def}.1).
  \item Then attach a copy of $E_2\cong \bar C_{j_2}$ to  $A\cupn E_1$ at a point of $A$ and at the point of $E_1$ lying over $c$. Also 
    attach a copy of $E_2\cong \bar C_{j_2}$ to  $B\cupn E_1$ at  2 points of $B$.
\end{enumerate}
After finitely many such steps, we achieve that $\# N^i_A(c, j_1, j_2)=\# N^i_B(c, j_1, j_2)$ for every $c, j_1, j_2$.

From now on we will not create more separating nodes.

Next we arrange that $\chi(\bar A^i_j)=\chi(\bar B^i_j)$ for every $j$.
Fix a curve $C_1$. If say $\chi(\bar A^i_j)>\chi(\bar B^i_j)$ for some $j\neq 1$ we use (\ref{elem.ch.over.c.def}.2) to attach a $\p^1$ at 2 points of $ \bar A^i_j$ and at
2 points of $ \bar B^i_1$.
This decreases $\chi(\bar A^i_j)$ by $1$ but leaves $\chi(\bar B^i_j)$ unaltered.
Repeating as necessary, and increasing the value of $i$,  we achieve that
$\chi(\bar A^i_j)=\chi(\bar B^i_j)$ for every $j\neq 1$.
Next note that
$$
\chi(A^i)=\tsum_j \chi(\bar A^i_j)-\tsum \# N^i_A(c, j_1, j_2),\qtq{and}
\chi(B^i)=\tsum_j \chi(\bar B^i_j)-\tsum \# N^i_B(c, j_1, j_2).
$$
We arranged that all but 1 of the terms in these identities match up, thus
$\chi(\bar A^i_1)=\chi(\bar B^i_1)$ also holds.

Next choose embeddings  $A^i\subset C\times \p^3$ and  $B^i\subset C\times \p^3$
that have the same  bidegree and such that 
$N^i_A(c, j_1, j_2)=N^i_B(c, j_1, j_2)$ for every $c, j_1, j_2$.

From now on we work separately with these 
$\bigl(\pi^i_A:\bar A^i_j\to \bar C_j,\pi^i_B:\bar B^i_j\to \bar C_j\bigr)$.
In (\ref{to.irred.curve.lem})  we prove a slightly stronger version of our claim, where all the steps leave the set of separated nodes fixed.
The 2-nodal deformation equivalence at the end fixes all the separated nodes.
If this holds, then the 2-nodal deformation equivalences over the individual 
$\bar C_j$ glue together to a 2-nodal deformation equivalence over $C$, as claimed. \qed

\begin{lem}  \label{to.irred.curve.lem}
  Let $C$ be a smooth, projective, irreducible  curve and
  $A^1, B^1\subset C\times \p^n$ two nodal subcurves
such that $\bigl(\bideg (A^1), \chi(A^1)\bigr)=\bigl(\bideg (B^1), \chi(B^1)\bigr)$.
Let $P\subset A^1\cap B^1$ be a finite set of points that are smooth on both $A^1$ and $B^1$.
  Then there is a sequence of elementary changes  (\ref{elem.ch.over.c.def}.1--3) fixing $P$ 
  such that  $\bigl(\pi^r_A:A^r\to C\bigr) =\bigl(\pi^r_B:B^r\to C\bigr)$.
\end{lem}

Proof.  First we achieve that  
$A^i\to C$ is finite and separable.
To do this, 
choose a smooth complete intersection surface   $A^1\subset S\subset C\times \p^3$ and $m>0$ such that  $\o_S(m)(-A^1)$ is very ample. A general section of it gives a smooth curve $F\subset S$.  We attach $F$ to $A^1$ along $F\cap A^1$ and to
$B^1$ along the same number if points as in (\ref{elem.ch.over.c.def}.1).
Now $A^2=A^1\cupn F$ is a very ample curve on $S$. It is thus linearly equivalent to a curve $A^3$ that is smooth, separable over $C$, and contains $P$.
This is (\ref{elem.ch.over.c.def}.3) with $E=\emptyset$.
We can now do the same for the $B$-curves.

Once $A^i, B^i$ are both finite and separable over $C$,
in characteristic 0 we could use \cite[1.1]{ha-ga-st}. It should work in any characteristic, but here is a weaker version that is good enough for us.

Choose a general projection
to $C\times \p^1$. Let the images be  $A', B'$, and  $P'\subset A'\cap B'$ the image of $P$.  We may assume that  $A', B'$ are smooth at $P'$.
Since  $A', B'$ have the same bidegree, they have 
 the same arithmetic genus. Thus we get the same number of new nodes $N'_A, N'_B$.
After adding the same number of vertical lines $L'_A=\cup_c \{c\}\times \p^1$, we may assume that  $A'+L'_A$ and $B'+L'_B$ are linearly equivalent, and  such that $|A'+L'_A|(-P')$ is very ample.
This gives the deformation equivalence of 
$A'+L'_A$ and of $B'+L'_B$. The lines  $L'_A, L'_B$ can be lifted to $L_A, L_B$
but we have a problem at the nodes $N'$ that are  smoothed in this deformation equivalence.
However, these can be lifted to the required deformation equivalence,
after a ramified double cover and adding lines that connect the preimages of the nodes as in (\ref{node.def.lift}). \qed

\begin{exmp} \label{node.def.lift}
  Consider the  smoothing of the node
  $(xy=t)\subset \a^3_{xyt}$. Take the ramified double cover  $t=s^2+s^3$.  (The cube is needed only in $\chr 2$.) Now we have   $(xy=s^2+s^3)\subset \a^3_{xys}$.
  Blowing up the origin, the central fiber is the birational transforms of the $x$ and $y$-axes, connected by the exceptional curve of the blow up.
\end{exmp}

\section{Gluing  of deformations}\label{glue.sec}

\begin{say}[Gluing deformations]\label{glue.def.say}
  As in (\ref{def.curves.defn}.1), consider 2 deformations
  $$
  B \longleftarrow S^j   \stackrel{\pi^j}{\longrightarrow}  B\times X
  \eqno{(\ref{glue.def.say}.1)}
  $$
    with sections $\sigma^j:B\to S^j$.

  If $\pi^1\circ \sigma^1\equiv \pi^2\circ \sigma^2$  then we can glue the 2
  deformations along the isomorphism
  $$
  S^1\supset \sigma^1(B)\cong B\cong \sigma^2(B)\subset S^2
  $$
  to get a new deformation
  $$
  B \longleftarrow  S^1\amalg_{\sigma} S^2   \stackrel{\pi}{\longrightarrow}  B\times X,
  \eqno{(\ref{glue.def.say}.2)}
  $$
  whose fiber over $b\in B$ is $S_b^1\amalg_{\sigma_b} S_b^2 $.
    Furthermore, if the (\ref{glue.def.say}.1) are 3-nodal,
  the $\sigma^j(b)$ is never a 3-node, and
  $\sigma^1(b), \sigma^1(b)$ are never both 2-nodes, then
  (\ref{glue.def.say}.2) is also 3-nodal.
\end{say}

The condition $\pi^1\circ \sigma^1\equiv \pi^2\circ \sigma^2$ is, however, very restrictive, so we are unlikely to be able to glue right away.
However, we can always connect any 2 points of $X$ by some auxiliary curve.
Our plan is to  glue both $S^j$ to the auxiliary family.

\begin{say}[Prelimnary steps]\label{prelsteps.say}{\ }

  {\it Creating sections \ref{prelsteps.say}.1.} 
  Given a deformation as in (\ref{def.curves.defn}.1) such that $S\to B$ has connected fibers. Let
  $p_i\in S_i$ be smooth points. Let $B'$ be a general member of a sufficiently  ample linear system   on $S$ passing through the $p_i$. After base change to
  $B'\to B$, the family $S'\to B'$ has a section connecting the points  
  $p'_i$.  If $S\to B$ is  3-nodal, then so is $S'\to B'$.
  (Note that $B'$ needs to avoid the 3-nodes, and that we pick up a new pivot points where $B'$ intersects the double curve of $S$.)

  If the $p_1, p_2$ are on different irreducible components of $S$, then $B'$ is reducible. However, if $B$ is  irreducible and $p_1, p_2$ are on the same irreducible component of $S$, then $B'$ is irreducible.

  Iterating this we can create any number of sections. We may assume that they meet only at smooth points and transversally. We can now blow up the intersection points to get a deformation with disjoint sections.

   {\it Creating isomorphic base curves \ref{prelsteps.say}.2.} 
  Given  pointed, connected curves  $(B^1, b^1_i)$ and
  $(B^2, b^2_i)$, let  $B$ be a general member of a sufficiently  ample linear system   on $B^1\times B^2$ passing through the $b_i:=(b^1_i, b^2_i)$.  The projections give base changes
  $\beta^1:(B, b_i)\to (B^1, b^1_i)$ and  $\beta^2:(B, b_i)\to (B^2, b^2_i)$.

  Thus, given 2 families  $S^j\to B^j$, after a base change we may assume that
  $B^1=B^2$. If the $S^j\to B^j$ are 3-nodal, so are the new families over $B$.
\end{say}

\begin{say}[Joining 2 points]
  Let $\tau_i:B\to \p^n$ be 2 morphisms form a curve $B$ to $\p^n$.
  For $b\in B$, let $L_b\subset \p^n$ be the line joining $\tau_1(b)$ and $\tau_2(b)$. This gives  $T:\p^1\times B\to \p^n$ whose restriction to the $0$ and $\infty$ sections give the $\tau_i$.
  There are, however, some degenerate cases.  First, if $\tau_1(b)=\tau_2(b)$ for every $b$.
  Second, if  $\tau_1(b)=\tau_2(b)$ holds only for finitely many $b_i$, we get
  the lines over $B\setminus\{b_i\}$ and then taking the limit should give the lines over the $b_i$. This works if the $b_i$ are smooth points, but not otherwise.

  We need a similar systematic way to connect 2 points in any variety $X$.
  Instead of lines, we use complete intersection curves. In order to avoid the 
  above problems with coinciding points, we use   complete intersection curves in $X\times \p^1$.
\end{say}

\begin{say}[$H$-complete intersection curves]\label{2.conn.fam.hci}
   Let $|H|$ be a sufficiently ample linear system on $X$ and $n:=\dim X$.
   (For example,  any $H\sim mA$  works, where $A$ is very ample and $m\geq 2$.)
   Consider $X\times \p^1$ with coordinate projections $p_i$.
   We have the linear systems $p_1^*|H|$ and $|H^*|:=p_1^*|H|+p_2^*|\o_{\p^1}(2)|$.

   We have a universal family
   $$
  M_H\stackrel{u}{\longleftarrow}  (T_H\subset U_H) \stackrel{\rho}{\longrightarrow}  X,
  \eqno{(\ref{2.conn.fam.hci}.1)}
    $$
    where  $M_H$ parametrizes objects 
    $
    (T\subset  L=H_1\cap\cdots\cap H_{n-1}\cap H^*_n),
    $ where $T$ is a length 2 subscheme of $X\times \p^1$, and
    $L$ is a smooth, complete intersection of 
    $H_1,\dots, H_{n-1}\in p_1^*|H|$  and of
    $H^*_n\in |H^*|$.
    We refer to this as the {\it \hci family.} Its members are double covers
    of $H$-complete intersection curves in $X$.

 Let $Y$ be a smooth, proper variety. The Hilbert scheme of length 2 subschemes with a marked point is $\hilb_{(2)}(Y)\cong B_{\Delta}(Y\times Y)$, the blow-up of the diagonal of $ Y\times Y$.
    
    For the \hci family,  $M_H$ is an open subset of a
    fiber bundle over
    $\hilb_{(2)}(X\times \p^1)$
    whose fiber is a product of projective spaces, and 
the induced map  $M_H\to \hilb_{(2)}(X\times \p^1)$ is surjective.
    
In particular, we have the following lifting property.

 \medskip

 {\it Claim \ref{2.conn.fam.hci}.2.} Let $B$ be a reduced curve,
 $b_j\in B$ a finite set of  points and
    $g:B\to \hilb_{(2)}(X\times \p^1)$ a morphism. For each  point
    $b_j$, choose $m_j\in M_H$ such that $\rho(T_{m_j})=g(b_j)$. Then there is a rational map $g':B\map M_H$ such that
    \begin{enumerate}
    \item[(a)]  $\rho\circ g'=g$,
    \item[(b)] $g'$ is defined at the  poins $b_j$, and $g'(b_j)=m_j$.
      \qed
      \end{enumerate}
    \end{say}

\begin{cor}\label{hci.2pt.cor}
  Let $X$ be a smooth, projective variety and fix an \hci family as in
  (\ref{2.conn.fam.hci}.1). Let $B$ be a reduced curve, $b_i\in B$ a finite set of  points  and
  $\tau^j:B\to X$ two morphisms.
  Choose liftings  $\tilde\tau^j:B\to X\times \p^1$ such that
  $\tilde\tau^1(b)=\tilde\tau^2(b)$ happens only at finitely many smooth points.
  Thus we have
  $g:=(\tilde\tau^1, \tilde\tau^2): B\to \hilb_{(2)}(X\times \p^1)$.
  For each $b_i$ choose \hci curves $L_i$ through $g(b_i)$. 
  
  Then there is a
  deformation 
  $$
  \begin{array}{ccccc}
   H(\tau)_i  & \subset  & H(\tau)  &  \stackrel{\pi}{\longrightarrow} & B\times X\\
  \downarrow &  &  \downarrow &  &  \downarrow \\
  b_i  & \in  & B   & = & B,
  \end{array}
  \eqno{(\ref{hci.2pt.cor}.1)}
  $$
  with 2 sections $\sigma^j:B\to H(\tau)$ such that
   \begin{enumerate}\setcounter{enumi}{1}
   \item  $\tau^j=\pi\circ \sigma^j$ for $j=1,2$,
   \item all but finitely many fibers of  $H(\tau) \to B$ are \hci curves, and
   \item $ H(\tau)_i =L_i$ for every $i$.
     \end{enumerate}
\end{cor}

Proof.  
The lifting of $g: B\to \hilb_{(2)}(X\times \p^1)$ to
$g':B\to  M_H$ is given by (\ref{2.conn.fam.hci}.2). The pull-back of the universal family by $g'$ gives $H(\tau)$. \qed

\medskip

We use \hci families to glue curves together. The situation is clear for disjoint curves.

\begin{defn}  \label{sliding.curves.defn}  Let $\pi_i: C_i\to X$ be reduced curves and
  $$
  {\mathbf L}:=\bigl\{(\tau_L: L\to X, p_1\neq p_2\in L)\}
  $$
  a set of 2-pointed curves. Let
  $\operatorname{Join}(C_1, r{\mathbf L}, C_2)$ denote the set of all curves of the form
  $$
  C_1\cup_{\sigma_1} (\amalg_{j\in J}L_{j})\cup_{\sigma_2}  C_2
  $$
  where   $\# J=r$, the $p_{ij}\in L$ are smooth, and 
  $\sigma_i: \amalg_j \{p_{ij}\}\into C_i^{\rm nodal}$ are injections
    such that $\tau_{L_j}(p_{ij})=\pi_i(\sigma_i(p_{ij}))$.
    Let $\operatorname{Join}^{\rm nodal}(C_1, r{\mathbf L}, C_2)$
    be the subset of those curves for which the $\sigma_i(p_{ij})\in C_i$ are smooth.
    \end{defn}

\begin{cor} \label{sliding.curves.lem}  Let $C_i\to X$ be  reduced curves and  $H$ an ample divisor class as in (\ref{2.conn.fam.hci}). Assume that the nodal loci of the $C_i$ are connected. Then
  \begin{enumerate}
  \item $\operatorname{Join}(C_1, rH^{\rm ci}, C_2)$ is connected, and
  \item $\operatorname{Join}^{\rm nodal}(C_1, rH^{\rm ci}, C_2)$ is a dense, open subset of it.  \qed
  \end{enumerate}
  \end{cor}

For families with some intersection points, we need to be  more careful. 

\begin{thm} \label{join.at.nodes.hci.thm}
  Let $X$ be a smooth, projective variety over an algebraically closed  field $K$, and 
  $\pi_i:(C_i\cupn D_i)\to X$ nodal curves mapping to $X$  for $i=1,2$.
   Assume that  the $C_i, D_i$ are connected, 
    $\#(C_1\cap D_1)=\#(C_2\cap D_2)$ (call this number $r$) and  we have
  3-nodal deformation equivalences
$$
  \begin{array}{ccccc}
   C_i & \subset  & S_C  &  \stackrel{\pi_C}{\longrightarrow} & B_C\times X\\
  \downarrow &  &  \downarrow &  &  \downarrow \\
  b_i  & \in  & B_C   & = & B_C,
  \end{array}
  %%\eqno{(\ref{def.curves.defn}.1)}
  \qtq{and}
  \begin{array}{ccccc}
   D_i & \subset  & S_D  &  \stackrel{\pi_D}{\longrightarrow} & B_D\times X\\
  \downarrow &  &  \downarrow &  &  \downarrow \\
  b_i  & \in  & B_D   & = & B_D.
  \end{array}
  $$
  Fix an \hci family of curves as in  (\ref{2.conn.fam.hci}.1).

  Then, after a base change  $\beta_C:B\to B_C$ and $\beta_D:B\to  B_D$,
  there is a 3-nodal deformation equivalence glued from
  $\beta_C^*(S_C)\to B$, $\beta_D^*(S_D)\to B$ and $r$ different \hci families.
\end{thm}

Proof.   By (\ref{prelsteps.say}), after a base change $\beta_C:B\to B_C$ and $\beta_D:B\to  B_D$ we may assume that there are 
disjoint  sections $\sigma^C_j:B\to S_C$ and  $\sigma^D_j:B\to S_D$
connecting the points in $C_1\cap D_1$ with the points in $C_2\cap D_2$.
This gives  $\tau^C_j:=\pi_C\circ \sigma^C_j:B\to  X$ and
$\tau^D_j:=\pi_D\circ \sigma^D_j:B\to  X$.

We apply (\ref{hci.2pt.cor})  to the pair of morphisms $\tau^C_j, \tau^D_j$.
We choose liftings $\tilde \tau^C_j:B\to  X\times\p^1$ and
$\tilde\tau^D_j:B\to  X\times\p^1$ such that
$\tilde \tau^C_j(b_1)=\tilde \tau^D_j(b_2)$ for every $j$.

Now use (\ref{hci.2pt.cor}) to create \hci families
$H^j\to B$ with  2 sections  $\sigma^{HC}_j:B\to H^j$ and
$\sigma^{HD}_j:B\to H^j$.

Now first glue $S_C$ to each $H^j$ using the $r$ section pairs $(\sigma^C_j, \sigma^{HC}_j)$.
Then glue the resulting surface to $S_D$ using the $r$ section pairs
$(\sigma^D_j, \sigma^{HD}_j)$. \qed

\begin{complem-p}\label{triple.glue.exmp}
  It may be worthwhile to write down explicitly what happens at the points
  $C_i\cap D_i$.

  By construction, at each point  we glue 3 smooth surfaces.  In suitable (formal or \'etale)
  local coordinates the gluing problem becomes
  $$
  \begin{array}{l}
  S_C\cong \a^2_{xt}\ni (0,t)\leftrightarrow  (g(t), t)\in \a^2_{zt}\cong H^j,\qtq{and}\\
  S_D\cong \a^2_{yt}\ni (0,t)\leftrightarrow  (h(t), t)\in \a^2_{zt}\cong H^j,
  \end{array}
  $$
  where $g(0)=h(0)=0$.
The $(t=0)$ fiber is the 3 coordinate axes, which can be given as 
 $$
 \left( \rank
  \left(
\begin{array}{ccc}
x & 0& z \\
0 & y& z 
\end{array}
\right)
\leq 1\right)\subset \a^3_{xyz}.
$$
By \cite[\S 5]{art-def}, any flat deformation is obtained  by varying the entries of the above matrix. 
 In our case,  the glued surface can be written as
 $$
 \left( \rank
  \left(
\begin{array}{ccc}
x  & 0& z{-}g(t)\\
0 & y & z{-}h(t) 
\end{array}
\right)
\leq 1\right)\subset \a^4_{xyzt}.
%%  \eqno{(\ref{triple.glue.exmp}.2)}
$$
\end{complem-p}

\section{Modifying deformation equivalences}\label{modify.sec}

\begin{thm}\label{tail.same.thm}
  Let $X$ be a smooth, projective variety over an algebraically closed  field $k$, and 
  $$
\begin{array}{ccccc}
C_i\cupn D_i & \subset  & S  &  \stackrel{\pi}{\longrightarrow} & B\times X\\
  \downarrow &  &  \downarrow &  &  \downarrow \\
  b_i  & \in  & B   & = & B
  \end{array}
$$
a 3-nodal deformation equivalence.
Assume that  the $C_i, D_i$ are connected, $\chi(D_1)=\chi(D_2)$ and 
 $\pi_*[D_1]=\pi_*[D_2]$. 

Then, there is a
3-nodal deformation equivalence
  $$
\begin{array}{ccccc}
C_i\cupn  R_i & \subset  & S_R  &  \stackrel{\pi_R}{\longrightarrow} & B_R\times X\\
  \downarrow &  &  \downarrow &  &  \downarrow \\
  b_i  & \in  & B_R   & = & B_R,
  \end{array}
$$
such that
\begin{enumerate}
\item $(\pi_R: R_1\to X)\cong (\pi_R: R_2\to X)$,
  \item $R_i\cap C_i=D_i\cap C_i$, and
  \item $S_R$ is glued from 4 types of pieces:
     pull-backs of $S$,
      deformation equivalences over   $\supp\pi(D_1)=\supp\pi(D_2)$, 
      trivial families  $C_i\times B'\to B'$, and
         \hci families.
   \end{enumerate}
\end{thm}

Proof.  We apply (\ref{to.curves.prop}) to
$$
\pi|_{D_1}:D_1\to \supp\pi(D_1)\qtq{and}
\pi|_{D_2}:D_2\to \supp\pi(D_2)=\supp\pi(D_1).
$$
Then we show that each of the steps (\ref{elem.ch.over.c.def}.1--4) can be lifted to
3-nodal deformation equivalences.

For (\ref{elem.ch.over.c.def}.1) this is achieved using
(\ref{join.at.nodes.hci.thm}): we glue $S\to B$ to the
trivial family $E\times B\to B$ using  additional \hci families.
For (\ref{elem.ch.over.c.def}.2) we glue
 $S\to B$ to the
trivial family $\p^1\times D_1\to D_1$,
again using  additional \hci families.

Steps (\ref{elem.ch.over.c.def}.3--4) are already deformation equivalences.

This is almost what we need, except that in (\ref{join.at.nodes.hci.thm}) we added a number of \hci families. Thus the fibers over $b_1$ and $b_2$ are of the form
$$
\bigl(C_i\cupn R_i\bigr)\amalg_{\sigma}\bigl(\cup_{j\in J}L_{ij}\bigr)
$$
where the  $L_{ij}$ are \hci curves  attached to nodes of $C_i\cupn R_i $ by $\sigma$. By construction, we have the same number of \hci curves for $i=1,2$.

We can now slide these \hci curves away from the nodes using
(\ref{sliding.curves.lem}),
to get the same  \hci curves for $i=1,2$. Thus fibers over $b_1$ and $b_2$  are now
$$
\bigl(C_i\cupn R_i\bigr)\cupn\bigl(\cup_{j\in J}L'_{j}\bigr),
$$
which is what we wanted.
\qed

  \begin{cor} \label{same.tail.lem}
    Let $X$ be a smooth, projective  variety over an algebraically closed  field $K$, and
    $\pi_i:A_i\to X$  connected, nodal curves  mapping to $X$. Assume that  $\chi(A_1)=\chi(A_2)$ and 
    $(\pi_1)_*[A_1]=(\pi_2)_*[A_2]$.
Then there is a  3-nodal deformation equivalence
  $$
\begin{array}{ccccc}
A_i\cupn L_i\cupn R_i & \subset  & S  &  \stackrel{\pi}{\longrightarrow} & B\times X\\
  \downarrow &  &  \downarrow &  &  \downarrow \\
  b_i  & \in  & B   & = & B,
\end{array}
\eqno{(\ref{same.tail.lem}.1)}
$$
such that
\begin{enumerate}\setcounter{enumi}{1}
\item $(\pi: R_1\to X)\cong (\pi: R_2\to X)$,
  \item $L_i$ is an \hci curve, meeting both $A_i$ and $R_i$ at a  single smooth  point, and 
  \item $S$ is  glued from 4 types of pieces as in 
(\ref{tail.same.thm}.3), (with $A_i$ replacing $D_i$).
  \end{enumerate}
\end{cor}

  Proof.
Pick any  \hci curve $L$ meeting both $A_1$ and $A_2$ at a  single smooth  point.
  There is a trivial deformation equivalence
  $(A_1\cupn L\cupn A_2)\times \p^1\to \p^1$. Apply (\ref{tail.same.thm}) with  $D_1=A_1$ and $D_2=A_2$.   \qed
\medskip

It is worth noting that  the fibers of $S\to B$ move very little in $\chow_1(X)$.
To state this, let 
$|H^{\rm ci}|\subset \chow_1(X)$ denote the closure of the locus of images of \hci curves. 

\begin{complem} \label{same.tail.lem.c}
  We can choose (\ref{same.tail.lem}.1) with the following property. There is an effective 1-cycle $Z$  and morphisms $\tau_i:B\to |H^{\rm ci}|$ such that
  $$
  \pi_*[S_b]=Z+\tsum_i L(\tau_i(b))
  \eqno{(\ref{same.tail.lem.c}.1)}
  $$
  for every $b\in B$.  (We stress that (\ref{same.tail.lem.c}.1) is an actual equality of cycles.) \qed
  \end{complem}

\section{Proof of Theorems~\ref{main.same.tail.thm}--\ref{main.same.tail.thm.k}}\label{main.pf.sec}

We start with a key special case, where we have a more precise answer.

  \begin{prop} \label{main.same.tail.prop}
    Let $X$ be a smooth, projective  variety over an algebraically closed  field $K$, and
    $\pi_i:C_i\to X$  connected, nodal curves  mapping to $X$. Assume that  $\chi(C_1)=\chi(C_2)$ and 
    $(\pi_1)_*[C_1]\sim_a(\pi_2)_*[C_2]$.
Then there is a  3-nodal deformation equivalence
  $$
\begin{array}{ccccc}
C_i\cupn L_i\cupn R_i & \subset  & S  &  \stackrel{\pi}{\longrightarrow} & B\times X\\
  \downarrow &  &  \downarrow &  &  \downarrow \\
  b_i  & \in  & B   & = & B,
\end{array}
\eqno{(\ref{main.same.tail.prop}.1)}
$$
such that
\begin{enumerate}\setcounter{enumi}{1}
  \item $(\pi:  L_1\cupn R_1\to X)\cong (\pi:  L_2\cupn R_2\to X)$, and
  \item $L_i$ is an \hci curve, meeting both $C_i$ and $R_i$ at a  single smooth  point.
\end{enumerate}
Furthermore, if the $C_i$ have a common irreducible component $D$, then
we can also achieve that
\begin{enumerate}\setcounter{enumi}{3}
\item $L_i$ meets $C_i$ at a point on $D$.
  \end{enumerate}
\end{prop}

  Proof.  Let  $ B\leftarrow T \to T\times X$ be a nodal deformation equivalence certifying
  $(\pi_1)_*[C_1]\sim_a(\pi_2)_*[C_2]$. We may assume that $T\to B$ has connected fibers. We also have the  trivial deformation equivalences
  $ C_i\times B\to B$. By (\ref{join.at.nodes.hci.thm}) we can glue these 3 together, using
  2 different \hci families to get a 3-nodal deformation equivalence
  $$
\begin{array}{ccccc}
T'_i & \subset  & T'  &  \stackrel{\pi'}{\longrightarrow} & B\times X\\
  \downarrow &  &  \downarrow &  &  \downarrow \\
  b_i  & \in  & B   & = & B,
\end{array}
\eqno{(\ref{main.same.tail.prop}.6)}
$$
such that  $T'_i=C_1\cupn L_1\cupn T_i \cupn L_2\cupn C_2$.

Apply (\ref{tail.same.thm}) with  $D_1= L_1\cupn T_i \cupn L_2\cupn C_2$ and
$D_2=C_1\cupn L_1\cupn T_i \cupn L_2$.   \qed
\medskip

\begin{say}[Proof of Theorem~\ref{main.same.tail.thm}]
  \label{main.same.tail.thm.p.23}  First apply 
(\ref{gc.chi.same.lem})  to get morphisms of nodal curves  $\bar\pi_i:C_i\cupn D\to X$ such that,
  the $C_i\cupn D $ are geometrically connected, and
   $\chi(C_i\cupn D)$ is independent of $i$.

  Now we are in the situation of (\ref{main.same.tail.prop}).
  Thus, for any two $i, j$  we have deformation equivalences between curves of the form
  $$
  (C_i\cupn D)\cupn L_{ij}\cupn R_{ij}\qtq{and}
  (C_j\cupn D)\cupn L_{ij}\cupn R_{ij}.
  $$
  We can choose the $L_{ij}$ to meet  the $C_i\cupn D$ and $C_j\cupn D$  at a point of $D$ by (\ref{main.same.tail.prop}.4).
  Finally, for any $\ell$, take
  $$
  (C_\ell\cupn D)\cupn \bigl(\amalg_{ij} L_{ij}\bigr)\cupn \bigl(\amalg_{ij}R_{ij}\bigr).
  \eqno{(\ref{main.same.tail.thm.p.23}.1)}
  $$
  We claim that they are 3-nodal deformation equivalent to each other.
  By symmetry consider $\ell=1,2$. Then we have
$$
  (C_\ell\cupn D)\cupn \bigl(\amalg_{ij} L_{ij}\bigr)\cupn \bigl(\amalg_{ij}R_{ij}\bigr)=
  R_{12}\cupn L_{12}\cupn (C_\ell\cupn D)\cupn \bigl(\amalg'_{ij} L_{ij}\bigr)\cupn \bigl(\amalg'_{ij}R_{ij}\bigr)
  \eqno{(\ref{main.same.tail.thm.p.23}.2)}
  $$
  where $\amalg'$ means that we omit the pair $ij=12$. Here
  $$
  R_{12}\cupn L_{12}\cupn (C_1\cupn D) \qtq{and}
  R_{12}\cupn L_{12}\cupn (C_2\cupn D)
  $$
  are 3-nodal deformation equivalent by construction.
  We can now use (\ref{sliding.curves.lem}) to add the
  $L_{ij}\cupn R_{ij}$ one at a time (for $ij\neq 12$).

  This gives a deformation equivalence that has some 3-nodal fibers, they can be eliminated using (\ref{jump.over.say}). \qed

  \end{say}

\begin{say}
[Proof of Theorem~\ref{main.same.tail.thm.k}]\label{main.same.tail.thm.p.4} 
We follow the arguments in (\ref{main.same.tail.thm.p.23}).  First apply 
(\ref{gc.chi.same.lem})  to get morphisms of nodal $k$-curves  $\bar\pi_i:C_i\cupn D\to X$ such that
the $C_i\cupn D $ are geometrically connected and
 $\chi(C_i\cupn D)$ is independent of $i$.

  We are in the situation of (\ref{main.same.tail.prop}--\ref{main.same.tail.thm.p.23}).
  Thus, working over the algebraic closure $K\supset k$,  we
  have the curves
  $$
  (C_\ell\cupn D)\cupn \bigl(\amalg_{ij} L_{ij}\bigr)\cupn \bigl(\amalg_{ij}R_{ij}\bigr).
  \eqno{(\ref{main.same.tail.thm.p.23}.1)}
  $$
  as in (\ref{main.same.tail.thm.p.23}.1).  Set $p_{ij}=D\cap L_{ij}$.

  Next fix a pair $ij$. The plan is to  replace $L_{ij}\cupn R_{ij}$ by all of its Galois conjugates   $(L_{ij}\cupn R_{ij})^{\sigma}$ and attach them to $D$ at the point $p_{ij}^{\sigma}$. In order to keep the new curve nodal, we need the various
  $p_{ij}^{\sigma}$ to be different from each other and from the points in $C_{\ell}\cap D$.

  Here it is important that we are free to move the curves  $L_{ij}$, and use
  (\ref{same.res.field.cor.1}) to guarantee a choice for which
  $L_{ij}\cupn R_{ij}$ is defined over $k(p_{ij})$, while avoiding any given finite subset of points of $D$.

  We can now do this for all the pairs $ij$ in some order. Of course we need to make sure that a new $p_{ij}$ is distinct form all nodes created earlier. \qed

  In this construction, the degrees $k(p_{ij})/k$ end up different.
  With a more careful choice in (\ref{same.res.field.cor.1}) we can achieve that these degrees are the same. It is, however, unlikely that one can choose the fields  $k(p_{ij})$ the same.
  \end{say}

 \begin{lem} \label{gc.chi.same.lem}
   Let $k$ be a perfect field and $L/k$ a finite field extension.
    Let $X$ be a smooth, projective, $k$-variety
  and  $\pi_i:C_i\to X_L$ (for $i\in I$) finitely many morphisms of reduced curves defined over $L$.
  Assume that either $\ind_{C_i}=1$  for every $i\in I$, or 
  $\chi(C_i)$ is independent of $i$. Then there is a
  smooth, geometrically  irreducible $k$-curve $\pi_D:D\to X$ and 
  morphisms of  curves  $\bar\pi_i:C_i\cupn D\to X_L$ such that
  \begin{enumerate}
  \item the $(C_i\cupn D)^{\rm nodal} $ are geometrically connected,
  \item $\chi(C_i\cupn D)$ is independent of $i$,  and
  \item $(\bar\pi_i|_D: D\to X)\cong (\pi_D: D\to X)$.
  \end{enumerate}
\end{lem}

 Note that the assumption on $\chi(C_i)$ is necessary by
 Example~\ref{quad.deg2.exmp}.

\medskip
Proof. Working in $Y:=X\times \p^3$, we may assume that the 
$C_i\subset Y$ are disjoint. Let $Z_i\subset C_i^{\rm sm}$ be  reduced subschemes that have at least 1 point on each irreducible component.
Let $D\subset Y$ be a smooth, complete intersection curve, defined over $k$, 
that contains the $Z_i$. Now attach $D$ to $C_i$ at $Z_i$
to get $C_i\cupn D$. Note that
$$
\chi(C_i\cupn D)=\chi(C_i)+\chi(D)-\deg Z_i.
$$
If  $\ind_{C_i}=1$  for every $i$ then we can choose the $Z_i$ such that
$\chi(C_i\cupn D)$ is independent of $i$.
We can always choose the $\deg Z_i$ to be the same. Then $\chi(C_i\cupn D)$ is independent of $i$, and the other conditions are clear. \qed

\section{Smoothing combs}\label{comb.sec}

Let $X$ be a smooth projective variety and $C\subset X$ a smooth curve with normal bundle $N_C$. If $N_C$ is generated by global sections and
$H^1(C, N_C)=0$, then $C$ can be deformed in all directions inside $X$, but otherwise $C$ may have no deformations at all.
The technique of combs, originated in \cite{KMM92a, ghs}, relies on the observation that if we attach enough smooth rational curves $R_i$ to $C$, and the
$N_{R_i}$ are generated by global sections, then the resulting  nodal curve $C\cupn \amalg_i R_i$ has a better deformation theory than $C$ itself. We call $C$ the {\it handle} and the $R_i$ the {\it teeth} of the {\it comb}  $C\cupn \amalg_i R_i$.
(See also \cite[Sec.6]{ar-ko} for an introductory treatment.)   The proof uses 2 properties of the curves $R$.  At a general point $c\in C$,
\begin{itemize}
\item the tangent space $T_cX$ is spanned by
  $T_cC$ and the $T_cR$ for all $R\ni c$, and
\item  $H^1(R, N_R(-c))=0$  for all $R\ni c$.
  \end{itemize}
If $X$ is src, then there are rational curves $R$ through every point and in every direction with ample $N_R$, as needed.
If $X\to B$ has smooth src general fibers, then we get the needed rational curves in the smooth fibers, but not in the singular fibers.
So the method only works for handles that dominate $B$.

Thus, in this case, we proceed in 2 steps. First we attach rational curves to \hci curves to achieve the needed properties. Then use these (reducible, high genus) curves as
teeth for handles that are contained in the singular fibers.
To accommodate both steps, we work out the details in a general setting.

  Each projective family of nodal curves $B_R\leftarrow S_R\to B\times X$
  can be realized as an embedded surface in $  B\times Y$ where $Y:=X\times \p^3$.
  We can thus view the deformation as  a geometrically connected curve
  $B_R\into \hilb_1(Y)$ that parametrizes nodal curves.

\begin{say}[Hilbert scheme]\label{hilb.defn} For a projective $k$-scheme $Y$, let
  $\hilb_1(Y)$ denote the Hilbert scheme of 1-dimensional, closed subschemes.
    Let $\hilb_1^{\rm nodal}(Y), \hilb_1^{\rm 3-nodal}(Y)$ denote the
  subschemes parametrizing nodal (resp.\ 3-nodal) curves.
   We use
  $\hilb_1^{\rm unobst}(Y)$ to denote the open 
   subscheme parametrizing reduced curves, all of  whose singularities have unobstructed deformations.

    For a curve $[C]\in \hilb_1^{\rm unobst}(Y)$,
  let $N_C:=\shom_C(I_C/I_C^{(2)}, \o_C)$ denote its {\it normal sheaf.}
  (It is locally free for nodes, but not for 3-nodes.) Then
  \begin{enumerate}
  \item   $H^0(C, N_C)$
  is the tangent space of $\hilb_1(Y)$ at $[C]$, and
  \item $H^1(C, N_C)$ is the obstruction space for deforming
  $C\subset Y$;
  \end{enumerate}
  see, for example, \cite[I.2.14]{rc-book}.

  We say that
  $C$ is {\it free} if $N_C$ is generated by global sections
  and  $H^1(C, N_C)=0$.
  Let  $C=C_1+ C_2$. Restriction gives a natural injection
  $I_C/I_C^{(2)}\into I_{C_1}/I_{C_1}^{(2)}+I_{C_2}/I_{C_2}^{(2)}$, hence we get
  $$
  0\to N_{C_1}+N_{C_2}\to N_{C_1+C_2}\to T\to 0,
  $$
  where $T$ is supported at $C_1\cap C_2$. In particular, we get the following.
  \medskip

  {\it Claim \ref{hilb.defn}.3.}  If $C_1, C_2$ are free and   $C_1+ C_2$ is locally unobstructed, then $C_1+ C_2$ is also free. \qed
  \medskip
  
  The open 
  subscheme parametrizing free curves is denoted by $\hilb_1^{\rm free}(Y)$.
  Note that $\hilb_1^{\rm free}(Y)$ is smooth.
It has the following general  properties.
  \medskip

  {\it Claim \ref{hilb.defn}.4.}  Every irreducible component of $\hilb_1^{\rm free}(Y)$ is a connected component. \qed
  
  \medskip

  {\it Claim \ref{hilb.defn}.5.} 
  Assume that $[Z_1]$ and $[Z_2]$ are smooth $k$-points on the same
  irreducible component of $\hilb_1(Y)$. Then $[Z_1]\sim_a [Z_2]$ over $k$. \qed

\end{say}

   \begin{say}[Families of teeth]  \label{curve.fam.defn}
    Let $X$ be a proper variety. A {\it family of pointed nodal curves} over $X$   is a diagram
  $$
  M\stackrel{\pi}{\longleftarrow}  C_M \stackrel{\tau}{\longrightarrow}  X
  \eqno{(\ref{curve.fam.defn}.1)}
    $$
    where $\pi:C_M\to M$ is a proper, flat morphism with nodal curve fibers, plus 
    a section $s:M\to C_M$ such  that $\pi$ is smooth along 
    $s(M)$.  

    For $p\in M$ let $\tau_p:C_p\to X$ denote the restriction of $\tau$ to the fiber over $p$. If every $\tau_p$ is an embedding, we have a flat family of nodal curves in $X$.

    For $x\in X$, the condition
    $\tau\circ s(p)=x$ defines $M_x\subset M$ and the
     the subfamily  $C_{M_x}\to M_x$.

       If $C_p\to X$ is an immersion at  $s(p)$ for
      every $p\in M$, then the derivative of $\tau_p$ at $s(p)$ gives 
      $d\tau: M\to \p_X(T_X)$.
      If $\tau_p: C_p\to X$ is an embedding, then we have the normal bundle
      $N_{C_p}$.

      The following properties are various versions of saying that the curves
       $\tau_p: C_p\to X$ `move freely' in $X$.
   \begin{enumerate}\setcounter{enumi}{1}
         \item For every $x\in X$, $\tau|_{M_x}: C_{M_x}\setminus s(M_x)\to X$ is equidimensional.
    \item     $d\tau: M\to \p_X(T_X)$ is  smooth and surjective.
    \item  $d\tau$ has geometrically irreducible fibers.
        \item $H^1\bigl(C_p, N_{C_p}(-s(p))\bigr)=0$ for every $p\in M$.
   \end{enumerate}
Our discussions in  (\ref{2.conn.fam.hci})  show that the family of pointed, smooth 
\hci curves satisfy (2--4), but they almost never satisfy (5).

However, if $X$ is separably rationally connected, then there is a way to create a family that also satisfies (5),

   \end{say}

\begin{say}[Combs]
  Let $X$ be a  proper variety
    over  an algebraically closed  field $K$. Let 
  ${\mathbf W}$ be a set of curves on $X$ and 
  $$
   M\stackrel{\pi}{\longleftarrow}  T_M \stackrel{\tau}{\longrightarrow}  X
   \qtq{with} s: M\to T_M
  $$
  a  family of  pointed curves as in (\ref{curve.fam.defn}.1); called teeth from now on, and denoted by ${\mathbf T}$. 
  A {\it comb} with handle in ${\mathbf W}$ and $r$ teeth in ${\mathbf T}$ is a curve
  $C\cupn(\amalg_{i\in I} T_{p_i})$, where $C\in {\mathbf W}$,
  $p_i\in M$ and 
  $\sigma: \amalg_i s(p_i)\into C^{\rm sm}$ is a gluing morphism such that
  $\tau(s(p_i))=\sigma(s(p_i))$ for every $i$.
  If we allow $\sigma: \amalg_i s(p_i)\into C^{\rm nodal}$, we get a
  {\it 3-nodal comb.}
  It is an  {\it embedded comb} in the natural map
  to $X$ is an embedding.
  Let $\ecomb({\mathbf W}, r{\mathbf T})$ denote the set of all such subschemes.

A comb is called {\it balanced} if each irreducible component of
  $C$ has at least $r/(2c)$ teeth attached at smooth points, where $c$ is the number of 
irreducible components of $C$. (2 is pretty arbitrary choice, any larger constant would work for us.)

These form a subset
  $\ecomb^{\rm bal}({\mathbf W},r{\mathbf T})\subset \ecomb({\mathbf W},r{\mathbf T})$.

  If ${\mathbf W}, {\mathbf T}$ are constructible subsets of $\hilb_1(X)$, then so are 
  $\ecomb({\mathbf W},r{\mathbf T})$ and $\ecomb^{\rm bal}({\mathbf W},r{\mathbf T})$.
\end{say}

The following two theorems are  slight generalizations of \cite[Sec.2]{ghs}, see also
\cite[Sec.6]{ar-ko}.

\begin{thm} \label{ghs.thm.1} Let $X$ be a smooth, projective variety over a perfect field $k$.
Assume that either $X$  is src and $\dim X\geq 3$, or  src in codimension 1 and 
$\dim X\geq 4$.
Let  ${\mathbf R}$ be an irreducible component of the space of almost very free morphisms  $\p^1\to X$.
   \begin{enumerate} %%\setcounter{enumi}{2}
  \item 
    Let ${\mathbf W}$ be any \hci family. Then there is an $r>0$ and a
    geometrically irreducible  open subset
   ${\mathbf T}\subset \ecomb({\mathbf W}, r{\mathbf R})$  that satisfies  (\ref{curve.fam.defn}.2--5).
  \end{enumerate}
  \end{thm}

Proof.  As we noted, the \hci family satisfies (\ref{curve.fam.defn}.2--3), which implies the same for
$\ecomb({\mathbf W}, r{\mathbf R})$. The key property 
(\ref{curve.fam.defn}.5) is proved in \cite[Sec.2]{ghs}.

One needs to pay  some attention to (\ref{curve.fam.defn}.4).
${\mathbf R}$ parametrizes certain maps $g:\p^1\to X$. Their lifts 
$\bar g:\p^1\to X\times \p^1$ are  parametrized by $\bar{\mathbf R}:={\mathbf R}\times \PGL_2$.
We have the universal rational curve map
$\rho:\p^1_{\bar{\mathbf R}}\to X\times \p^1$. Given an \hci  curve  $L\subset X\times \p^1$, the curves $R\cupn L$ are parametrized by
$\rho^{-1}(L)$. By Bertini's connectedness theorem, $\rho^{-1}(L)$
is geometrically irreducible for general $L$; see \cite[4.10 and 6.10]{MR725671}.

If ${\mathbf R}$ itself is  geometrically irreducible, then
$\ecomb({\mathbf W}, r{\mathbf R})$ is  geometrically irreducible, and we are done. (In fact, one can always choose such an  ${\mathbf R}$.)
In general, if ${\mathbf R}$ decomposes into geometrically irreducible
components ${\mathbf R}_i$, then ${\mathbf T}$ consist of those combs that have the same number of teeth from each  ${\mathbf R}_i$. 
\qed

\medskip

Next we use ${\mathbf T}$ obtained in (\ref{ghs.thm.1}) as teeth of the next comb construction, but first a definition.

\begin{defn}\label{dfree.comb.defn}
  A comb $C\cupn(\amalg_i T_i)\subset X$ is {\it $d$-free along the handle}
  if for every subcurve $C_1\subset C$, 
  $$
  H^1\bigl(C_1\cupn(\amalg_i T_i), N_{C_1\cupn(\amalg_i T_i)}(-D)\bigr)=0
  \eqno{(\ref{dfree.comb.defn}.1)}
  $$
  for every  effective Cartier divisor $D\subset C_1^{\rm sm}$ of degree $\leq d$.  

  Given ${\mathbf W}, {\mathbf T}$, let 
$\ecomb^{\rm d-free}({\mathbf W}, r{\mathbf T})$ denote the set of all balanced, 3-nodal  combs that are free and  $d$-free along their handle.

  {\it Comments.} Requiring (\ref{dfree.comb.defn}.1) for all  $C_1\subset C$ will be important in
(\ref{main.k.src.thm.cor.2.pf}). 

  It would have been simpler to  require (\ref{dfree.comb.defn}.1)  for every
  $D\subset \bigl(C\cupn(\amalg T_i)\bigr)^{\rm sm}$. However,
  the teeth $T_i$ that we use contain the rational curves $R_i$
  of (\ref{ghs.thm.1}) as irreducible components. If $X$ is src, we can choose
  these $R_i$ to have arbitrarily positive normal bundle. If $X$ is only
  src in codimension 1, then the $R_i$ are contained in the fibers of  $X\to B$, thus their normal bundle   always has a trivial summand. So
  $H^1\bigl(R_i, N_{R_i}(-2)\bigr)$ is never 0.
  We could fix this by replacing the reducible curves 
  ${\mathbf T}$ obtained in (\ref{ghs.thm.1}) by their general smooth deformations. This would necessitate some extra steps to work over nonclosed fields.
  \end{defn}

\begin{thm} \label{ghs.thm} Let $Y$ be a smooth, projective variety and
  ${\mathbf W}\subset \hilb_1^{\rm unobst}(Y)$ a bounded set.
  Let ${\mathbf T}$ be a family of teeth satisfying (\ref{curve.fam.defn}.2--5) and fix $d>0$.
  
  Then  there is an $r>0$ such that, for every $C\in {\mathbf W}$,
  $\ecomb^{\rm d-free}(C, r{\mathbf T})$ is open and dense in 
  $\ecomb^{\rm bal}(C,r{\mathbf T})$.
  \end{thm}

Proof.  Choose  $C\in {\mathbf W}$ and let $D\subset C^{\rm sm}$ be an effective Cartier divisor of degree $\leq d$. 
$H^1\bigl(C, N_C(-D)\bigr)$ is dual to
$\Hom_C\bigl(N_C(-D),\omega_C\bigr)$.

The computation of \cite[Sec.2]{ghs}, which is done at smooth points of $C$, shows that if $\eta\in \Hom_C\bigl(N_C(-D),\omega_C\bigr)$, then attaching a
 tooth $T_c$ at a  point of $c\in \supp(\eta)$ kills $\eta$, provided $\eta$ is nonzero on the tangent space of $T_c$ at $c$.  Thus there is an
$r_0$ such that if we attach at least $r_0$ general teeth at  general points
to all irreducible components on $C$, then $H^1\bigl(C, N_C(-D)\bigr)$ gets killed.
Now take  $r\geq 2cr_0$,  where $c$ is the number of 
irreducible components of $C$.

Finally note that both the dimension of $H^1\bigl(C, N_C(-D)\bigr)$ and the number of 
irreducible components are uniformly bounded in a  bounded set of curves. 
\qed

\medskip

The main technical result is the following.

\begin{thm}\label{main.src.comb.thm}
  Let $Y$ be a smooth, projective variety of dimension $\geq 3$ over  a perfect field $k$.
  Let ${\mathbf W}\subset \hilb_1^{\rm nodal}(Y)$ be a   bounded, locally closed subset, and ${\mathbf T}$  a family of pointed curves satisfying
  (\ref{curve.fam.defn}.2--5). Fix $d>0$. 

  Then there is an $r>0$ such that the forgetful map
  $$
  \Pi: \ecomb^{\rm d-free}({\mathbf W}, r{\mathbf T})\longrightarrow
     {\mathbf W}
     $$
   \begin{enumerate}
   \item  has connected fibers,
     \item satisfies the curve lifting property (\ref{c.lift.prop.say}), and
     \item surjective on $L$-points for every $k\subset L\subset K$. 
   \end{enumerate}
\end{thm}

\begin{cor}\label{main.src.comb.thm.c}Using the notation of
  (\ref{main.src.comb.thm}), assume in addition that 
  ${\mathbf W}$ is geometrically connected.
  Then $\ecomb^{\rm d-free}({\mathbf W}, r{\mathbf T})$ is 
  contained in a single  geometrically irreducible component    $\env({\mathbf W},r{\mathbf T})$ of $\hilb_1(Y)$; called its  {\em envelope.}
\end{cor}

Proof. $\ecomb^{\rm d-free}({\mathbf W}, r{\mathbf T})$ is
geometrically connected by (\ref{c.lift.prop.say}).
Thus the unique  connected component of $\hilb_1^{\rm free}(Y)$
containing it is  geometrically irreducible by (\ref{hilb.defn}.4). \qed

\begin{cor}\label{main.src.comb.thm.c.2} Using the notation of
  (\ref{main.src.comb.thm}), assume in addition that 
  ${\mathbf W}$ is geometrically connected.
  If $[C_1], [C_2]\in {\mathbf W}(L)$, then
  $C_1$ and $C_2$ are algebraically equivalent over $L$.
\end{cor}

Proof. By (\ref{main.src.comb.thm}.3) there are free combs
$C_i\cupn (\amalg_j T_{ij})$ defined over $L$. By
(\ref{hilb.defn}.5), the $C_i\cupn (\amalg_j T_{ij})$ are
algebraically equivalent over $L$. The $\amalg_j T_{ij}$
are also smooth $L$-points of the $r$th symmetric power of
${\mathbf T}$, hence the $\amalg_j T_{ij}$ are also
algebraically equivalent to each other over $L$. \qed

\begin{say}[Proof of \ref{main.src.comb.thm}]\label{main.src.comb.thm.pf}
  Let $C\in {\mathbf W}$ with irreducible components $\{C_i: i\in I\}$.
  The combs where we attach $m_i$ teeth to $C_i$ form an irreducible family.
  By (\ref{ghs.thm}) there is an $r_0$ such that general such combs
  are $d$-free if $m_i\geq r_0$ for every $i$. If $m_i>r_0$ for every $i$, then,
  as we slide one of the teeth across a node of $C$, we get a 3-nodal comb, which is still $d$-free be  (\ref{hilb.defn}.3). Thus, if $r\geq (r_0+1)\cdot\#I$,
then $\ecomb^{\rm d-free}(C, r{\mathbf T})  $ is connected. 

In order to check the curve lifting property (\ref{c.lift.prop.say})
 let $(B,b_1, b_2)$ be an irreducible, smooth, 2-pointed  curve, and
$\tau:(B,b_1, b_2) \to {\mathbf W}$ a morphism. By pull-back we get a family  $C_B\to B$.
 Attach $m$ general teeth to each irreducible component of the $C_i$ at the points $c_{ij}\in C_i$.  After a further base change, we may assume to have sections  $\sigma_{ij}:B'\to C'_B$, such  that
 $\sigma_{ij}(b'_i)=c_{ij}$.   We keep $B'$ irreducible, thus we can not control
 on which irreducible components of $C_{3-i}$ the points $\sigma_{3-i,j}(b_i)$ are on. However, after attaching further teeth, we have a family in
 $\ecomb^{\rm d-free}(C, r{\mathbf T})  $.

 Assume next that $C$ is defined over a field $L$ and
 $C\cupn (\amalg_j T_j)$ is $d$-free and defined over $K$.
 Adding all the $\gal(K/L)$-conjugate teeth, we get a
 comb $C\cupn (\amalg_{j,\sigma} T_j^\sigma)$ which
 is $d$-free and defined over $L$. We need to make sure that the $T_j^\sigma$ are attached to $C$ at distinct points and that their number is uniformly bounded.
Such a choice of the $T_j$ is possible by (\ref{same.res.f.lem}).
 \qed

 \end{say}

\begin{say}[Curve lifting property]\label{c.lift.prop.say}
  Let $K$ be algebraically closed and 
   $g:Y\to X$  a morphism of $K$-schemes of finite type.  
  If $g$ is not proper, it can happen that $X$ is  connected,
the fibers of $g$ are nonempty and connected,
  but $Y$ is not. The implication, however, holds, if $g$ has the following property.

  Let $(B, b_1, b_2)$ be a 2-pointed, irreducible (not necessarily proper) curve and
  $\tau: (B, b_1, b_2)\to X$ a morphism. Then there is
  2-pointed, irreducible curve $(B', b'_1, b'_2)$ and a commutative diagram
  $$
  \begin{array}{ccc}
    (B', b'_1, b'_2)  &  \stackrel{\tau'}{\to} & Y\\
    \downarrow &&   \hphantom{g}\downarrow g\\
    (B, b_1, b_2)  &  \stackrel{\tau}{\to} & X.
  \end{array}
  $$
  \end{say}

\section{Arithmetic applications}\label{arith.sec}

First we discuss the index and the genera of curves.

\begin{say}[Euler characteristic of 1-cycles]\label{ind.ech.defn}
  Let $k$ be a  field, $X$ a proper $k$-scheme and
  $\ind_X$ the {\it index} of $X$, that is,  the gcd of the  degrees of all 0-cycles on $X$.

  Let $C\subset X$ be an irreducible curve.
  We check in (\ref{ind.ech.say}) that $\ind_X\mid 2\chi(C)$ and set
  $$
  \chi^*(C):= \tfrac{2\chi(C)}{\ind_X}\mod 2.
  $$
  Thus, if  $ \ind_X$ is odd then $\chi^*(C)=0$ for every $C$.
  
  We extend $\chi^*$  to  reducible curves and 1-cycles  $Z=\tsum_i d_iC_i$ by linearity.

 We check in (\ref{chi.alg.eq.lem}) that if  $Z$ is algebraically equivalent to 0 then  $\chi^*(Z)=0$. Thus $\chi^*$ descends to a homomorphism
 $\chi^*: A_1(X)\to \z/2\z $.
\end{say}

 \begin{say}\label{ind.ech.say}
  Let $k$ be a  field, $X$ a proper $k$-scheme and
  $\ind_X$ the  index of $X$.

  Let $p:C\to  X$ be a proper, reduced curve mapping to $X$ and $\pi:\bar C\to C$  the normalization. Then $p_*(\pi_*\o_{\bar C}/\o_C)$ is a 
0-cycle on $X$, hence its degree is divisible by $\ind_X$. Thus
$\chi(\bar C)\equiv \chi(C)\mod \ind_X$.

The canonical class of $\bar C$ is represented by a 0-cycle, and
$\deg K_{\bar C}=-2\chi(\bar C)$. Thus $\ind_X\mid 2\chi(\bar C)$.
We have 2 possibilities:
  \begin{enumerate}
  \item if $\ind_X$ is odd then $\chi(C)\equiv 0 \mod \ind_X$, and
  \item if $\ind_X$ is even then $\chi(C)\equiv \chi^*(C)(\ind_X/2)\mod \ind_X$.   \end{enumerate}  
 (This is related to the  elw indices
  defined in \cite{k-elw}. For $d$-cycles, we get a similar invariant
  modulo  $\gcd\bigl(\operatorname{elw}_{d-1}, \mu(\operatorname{Td}_d)\bigr)$
  where $\mu(\operatorname{Td}_d)$ is the denominator appearing in the Todd class in dimension $d$.)
 \end{say}
   
 \begin{prop}\label{chi.alg.eq.lem}  (cf.\ \cite[3--5]{k-elw})
   Let  $X$ be a proper $k$-scheme.
   Then the function on 1-cycles
   $Z\mapsto \chi^*(Z)\in \z/2\z $ 
  is preserved by algebraic equivalence.
\end{prop}

Proof. Let $p:C\to X$ be a proper  morphism from a nonsingular curve to $X$
such that
$p_*[C]=[Z]$. The Euler characteristic of fibers of flat morphisms is  a locally constant function on the base,  so it is enough to prove that
$$
\chi(C)\equiv \tsum_i d_i\chi(Z_i)\mod \ind_X.
\eqno{(\ref{chi.alg.eq.lem}.1)}
$$
By linearity, it is enough to prove this for a morphism  
  $\tau: C\to B$  of nonsingular, irreducible curves  over $X$. 

If $\tau(C)$ has dimension 0, then
$\chi(C)=\deg \tau_*\o_C-\deg R^1\tau_*\o_C$, and both terms on the right hand side are divisible by $\ind_X$.
If $\tau$ is finite then
$$
\chi(C)=\chi(\tau_*\o_{C})=\deg \tau_*\o_{C}+\deg \tau \cdot \chi(\o_B).
$$
Here $\deg \tau_*\o_{C}$ is  divisible by $\ind_X$. \qed

\begin{lem} \label{eq.chi.lem}
  Let $X$ be a smooth, projective variety over a perfect field $k$ with 1-cycles $Z_1, Z_2$ such that $\chi^*(Z_1)=\chi^*(Z_2)$.
  Let $A\subset X$ be a smooth, curve meeting all irreducible components of $Z_1\cup Z_2$ at smooth points and such that $\ind_X=\ind_A$.

  Then there are connected, nodal $k$-curves $\pi_i:C_i\to X$ such that
  $(\pi_i)_*[C_i]=Z_i+2A$ and $\chi(C_1)=\chi(C_2)$.
\end{lem}

Proof.  We have  connected, nodal $k$-curves $\rho_i:D_i\to X$ 
such that
$(\rho_i)_*[D_i]=\red (Z_i)+A$. If $D_{ij}\subset D_i$ is an irreducible component, and $(\rho_i)_*[D_{ij}]$ appears in $Z_i$  with multiplicity $d_{ij}$, then we replace  $D_{ij}$ with a degree  $d_{ij}$ cyclic cover  $C_{ij}\to  D_{ij}$  that ramifies along all the nodes. Now we have  connected, nodal $k$-curves $C'_i=\cup C_{ij}$ such that
$(\pi'_i)_*[C'_i]=Z_i$.

Note that $\ind_X$ divides 
$\chi(C'_1)-\chi(C'_2)$. Also, since $\ind_X=\ind_A$, the curve
$A$ has double covers  $A_i\to A$ such that
$\chi(A_2)-\chi(A_1)=\chi(C'_1)-\chi(C'_2)$.
We now get $C_i$ by gluing $C'_i$ to $A_i$. \qed

\medskip

The stronger form of  Theorem~\ref{main.thm.1} is the following.

\begin{thm}\label{main.thm.1.str}
  Let $X_k$ be a smooth, projective   variety over  a perfect field $k$ with algebraic closure $K$. Assume that $X_K$ is  src in codimension 1.
  Then 
  $$
  A_1(X_k)\to A_1(X_K)+\z/2\z\qtq{given by}
  [Z_k]\mapsto  \bigl([Z_K], \chi^*(Z_k)\bigr)
  $$
  is injective.
\end{thm}

Proof.  If $Z_k\sim_a 0$ then also  $Z_K\sim_a 0$, and
$\chi^*(Z)=0$ by (\ref{chi.alg.eq.lem}).

Conversely, assume that $\chi^*(Z)=0$. By (\ref{eq.chi.lem})
there are connected, nodal $k$-curves $\pi_i:C_i\to X$ such that
$Z\sim_a (\pi_1)_*[C_1]-(\pi_2)_*[C_2]$ and $\chi(C_1)=\chi(C_2)$.

Then, by  Theorem~\ref{main.same.tail.thm.k}, there is a nodal deformation equivalence
  $$
\begin{array}{ccccc}
C_i\cupn  R & \subset  & S  &  \stackrel{\pi}{\longrightarrow} & B\times X\\
  \downarrow &  &  \downarrow &  &  \downarrow \\
  b_i  & \in  & B   & = & B,
\end{array}
%%\eqno{(\ref{main.same.tail.thm}.1)}
$$
where the $\pi|_{C_i\cupn R}:C_i\cupn R\to X$ are  defined over $k$.
Note that 
$$
(\pi_1)_*[C_1\cupn  R]-(\pi_2)_*[C_2\cupn  R]=
(\pi_1)_*[C_1]-(\pi_2)_*[C_2]\sim_a Z.
$$
Finally  $(\pi_1)_*[C_1\cupn  R]-(\pi_2)_*[C_2\cupn  R]\sim_a 0$ by
(\ref{main.src.comb.thm.c.2}) applied to $C_i\cupn R$. \qed

\begin{say}[Proof of Theorem~\ref{main.k.src.thm.cor.3}]
 Injectivity follows from (\ref{main.thm.1}.1). 

 Let $Z$ be a 1-cycle on $X_K$ such that $Z\sim_a Z^\sigma$
 for every $\sigma\in \gal(K/k)$. By Theorem~\ref{main.same.tail.thm.k} we may certify $Z\sim_a Z^\sigma$  by a nodal curve  $C\cupn R\subset X\times \p^3$ that  is  nodal 
 deformation equivalent to  $C^\sigma\cupn R$ for every $\sigma\in \gal(K/k)$.
 The union of these nodal equivalences and their $\gal(K/k)$-conjugates
 gives a geometrically connected ${\mathbf W}\subset \hilb_1( X\times \p^3)$
 which contains all the $C^\sigma\cupn R$.

 Choose a geometrically irreducible family of teeth ${\mathbf T}$, and apply
 (\ref{main.src.comb.thm.c}). We get the envelope
 $\env({\mathbf W},r{\mathbf T})$ that contains a free comb
 $(C\cupn R)\cupn(\amalg_j T_j)$  and all of its
 $\gal(K/k)$-conjugates  $\bigl((C\cupn R)\cupn(\amalg_j T_j)\bigr)^\sigma$.
 Thus $\env({\mathbf W},r{\mathbf T})$ is geometrically irreducible.

By assumption there are 0-cycles 
$P$ on ${\mathbf T}$ and  $Q$ on $\env({\mathbf W},r{\mathbf T})$ of degree 1. Let
$T_P$ and $E_Q$ be the corresponding 1-cycles on $X\times \p^3$.
Then $[E_Q]-r[T_P]$ is a 1-cycle on  $X\times \p^3$ that is
algebraically equivalent to $C\cupn R$ over $K$. Thus
$[E_Q]-r[T_P]-[R]$ is
algebraically equivalent to $C$ over $K$.
Projecting it to $X$ gives a 1-cycle on  $X$ that is
algebraically equivalent to $[Z]$ over $K$. \qed

\end{say}

\begin{say}[Proof of Theorem~\ref{main.k.src.thm.cor.2}]
  \label{main.k.src.thm.cor.2.pf}
Pick $[Z_0]\in  A_1(X_0)$.
We can represent it as the projection of a nodal curve $C_0\subset Y:=X\times \p^3$.
By (\ref{ghs.thm}) there is a free comb  $C_0\cupn(\amalg_j T_j)$ where
$\amalg_j T_j$ is defined over $k$. 
Thus
$\hilb_1(Y/R)\to\spec R$ is smooth at 
$[C_0\cupn(\amalg_j T_j)]$ and also at $[\amalg_j T_j]$. Since $R$ is Henselian, these lift to
$[E_g], [T_g]\in \hilb_1(Y_g)$. Thus $[E_g]- [T_g]$ shows that $A_1(X_g)\to A_1(X_0)$ is surjective.

For injectivity, pick $Z^1_g, Z^2_g$ with specializations
$Z^1_0\sim_a Z^2_0$. Then  $\ind_{X_0}\mid \chi(Z^1_0)-\chi(Z^2_0)$ by
(\ref{main.thm.1}.2). 
Note that $\ind_{X_0}=\ind_{X_g}$, so 
$\ind_{X_g}\mid \chi(Z^1_g)-\chi(Z^2_g)$. Thus,  by (\ref{main.thm.1}.2), it is enough to show that $Z^1_g$ and $Z^2_g$ become algebraically equivalent after a field extension of $k(g)$.

By (\ref{gc.chi.same.lem}) we can choose nodal representatives
$C^i_g\to Z^i_g$ such that $\chi(C^1_g)=\chi(C^2_g)$. 
After semistable reduction, we may replace the $C^i_g$ by flat families of nodal curves $C^i_R\subset X_R\times_R \p^3_R$.
We apply Theorem~\ref{main.same.tail.thm} to get
$$
\begin{array}{ccccc}
C^i_0\cupn  E & \subset  & S  &  \stackrel{}{\hookrightarrow} & B\times X_0\times \p^3\\
  \downarrow &  &  \downarrow &  &  \downarrow \\
  b_i  & \in  & B   & = & B.
\end{array}
%%\eqno{(\ref{main.same.tail.thm}.1)}
$$

Now apply (\ref{main.src.comb.thm.c})
to ${\mathbf W}:=B$ with $d:=\max\{\# (C^1_0\cap E), \# (C^2_0\cap E)\}$.
We get   $d$-free combs
$(C^i_0\cupn E)\cupn (\amalg_j T^i_j)$ that  are in the same irreducible component of $\hilb_1^{\rm free}(Y_0)$.

The $C^i_0\cap E$ are smooth points on $C^i_0$, so can be lifted
to a union of sections  $\Sigma^i\subset C^i_R$. 
We can also view  $C^i_0\cap E$, as a Cartier divisor on $E$.
Since the comb is $d$-free, 
$$
H^1\bigl(E\cupn (\amalg_j T^i_j), N_{E\cupn (\amalg_j T^i_j)}(-(C^i_0\cap E))\bigr)=0,
$$
hence every deformation of  $C^i_0\cap E$ extends to a deformation of
$E\cupn (\amalg_j T^i_j)$.  In particular, 
each $E\cupn (\amalg_j T^i_j)$ can be lifted to  a flat family of curves $V^i_R$
containing $\Sigma^i$.  Then 
  $C^i_R\cup V^i_R$ is  flat over $R$.

Since the $E\cupn (\amalg_j T^i_j)=V^i_0$ are
smooth points on same irreducible component of $\hilb_1^{\rm free}(Y_0)$,
the same holds for the  general fibers $V^i_g$.
Hence $V^1_g\sim_a  V^2_g  $.

The same argument shows that
$C^1_g\cup V^1_g\sim_a C^2_g\cup V^2_g$. Therefore
$C^1_g\sim_a C^2_g$. \qed
\end{say}

\section{Loops in $\chow_1$}\label{chow.loop.sec}

\begin{defn}  Let $X$ be a  scheme. A length $m$ {\it chain}
   in $X$ is a
  collection of  morphisms
  $$
 {\mathcal L}:=\bigl\{ \pi^i:(W^i, Z^i_1, Z^i_2)\to X \colon i=1,\dots, m\bigr\},
  $$
  where the $W^i$ are irreducible schemes,
  $Z^i_1, Z^i_2$ are irreducible, closed subsets and
  $\pi^i(Z^i_2)=\pi^{i+1}(Z^{i+1}_1)$  (as sets).

  ${\mathcal L}$ is a   {\it loop}  if $\pi^m(Z^m_2)=\pi^{1}(Z^{1}_1)$ also holds. A loop is {\it contractible} if all the $\pi^i(W^i)$ are contained in the same irreducible component of $X$.

  A connected, pointed  curve
 $\pi_B^i:(B^i, p^i_1, p^i_2) \to (W^i, Z^i_1, Z^i_2)$ is a  {\it representative.}  Assume that  we also have 
  connected, pointed  curves
  $\pi_C^i:(C^i, q^i_1, q^i_2) \to p^i(Z^i_2)$ such
 that $\pi^i\circ \pi_B^i(p^i_2)=\pi_C^i(q^i_1)$ and
$\pi^{i+1}\circ \pi_B^{i+1}(p^i_1)=\pi_C^i(q^i_2)$ (indexing modulo $m$). 
 Alternating the $B^i$ with the $C^i$ gives a loop in $X$. We call it 
 a {\it representative} of ${\mathcal L}$.
 \end{defn}

\begin{defn}  Let $X$ be a proper scheme over a field and
  $$
 {\mathcal L}:=\bigl\{ \pi^i:(W^i, Z^i_1, Z^i_2)\to X \colon i=1,\dots, m\bigr\},
  $$
a loop in $\chow_1(X)$. Let $Z_V\subset V\subset \chow_1(X)$ be irreducible subvarieties. The {\it translation} of  ${\mathcal L}$
by $ (V\supset Z_V)$ is the loop made from the 
$$
(W^i\oplus V, Z^i_1\oplus Z_V, Z^i_2\oplus Z_V) \to \chow_1(X),
$$
where for $U_1, U_2\subset \chow_1(X)$, we set    $U_1\oplus U_2:=\{Z_1+Z_2: Z_i\in U_i\}$. 
\end{defn}

  \begin{thm} Let $X$ be a smooth, projective variety, and
    ${\mathcal L}$ a loop in $\chow_1(X)$.
    Then there is a 1-cycle $Z$ 
    and $r>0$ such that the translation of ${\mathcal L}$ by
      $[Z]\oplus r|H^{\rm ci}|$ is representable by
      a loop in  $\hilb_1^{\rm nodal}(X\times \p^3)$.
  \end{thm}

  Proof.
Choose points $z^i_j\in Z^i_j$. By (\ref{certifies.defn}), there are  2-nodal deformation equivalences
  $$
  (B^i, b^i_1, b^i_2)\stackrel{}{\longleftarrow} S^i
  \stackrel{\pi^i}{\longrightarrow}  B^i\times X
  $$
  representing $(W^i, z^i_1, z^i_2)$.
  After adding a suitable 1-cycle $D$ and using (\ref{gc.chi.same.lem}),
  we may assume that the $S^i\to B^i$ have connected fibers and their
  Euler characteristic is independent of $i$.

  The problem is that, while $S^i_2\to X$ and $S^{i+1}_1\to X$ give the same 1-cycles on $X$, they are not isomorphic as maps.

  Next we apply 
  (\ref{same.tail.lem}--\ref{same.tail.lem.c})
  to each $\pi^i:S^i_2\to X$ and $\pi^{i+1}:S^{i+1}_1\to X$.
  We get nodal curves $R_i$ and 3-nodal deformation equivalences 
  $$
  (C^i, c^i_1, c^i_2)\stackrel{}{\longleftarrow} T^i
  \stackrel{\tau^i}{\longrightarrow}  C^i\times X
  $$
  such that
  $$
  T^i_1\cong S^i_2\cupn L_i\cupn R_i \qtq{and} T^i_2\cong S^{i+1}_1\cupn L_i\cupn R_i.
  $$
  We can now assemble a 3-nodal  representative loop from the pieces
  $$
  S^i\cupn_j \bigl(L_j\cupn R_j)
  \qtq{and}
  T^i\cupn_{j\neq i} \bigl(L_j\cupn R_j). 
  $$
  Now take $[Z]=[D]+\sum_j [\tau^j_*(R_j)]$.
  We can turn it nodal using  (\ref{jump.over.say}). \qed

  \medskip

  Combining this with (\ref{main.src.comb.thm.c}) we get the following.

  \begin{cor} Let $X$ be a smooth, projective  variety, and
    ${\mathcal L}$ a loop in $\chow_1(X)$.
    Assume that $X$ is   src in codimension 1.

    Then a suitable translate of  ${\mathcal L}$ is
    contractible. \qed
  \end{cor}

\section{Controlling residue fields}\label{res.f.sec}

\medskip

In Sections~\ref{main.pf.sec}--\ref{comb.sec} we need the following to control residue fields.

\begin{lem}\label{same.res.f.lem} Let $k$ be  a perfect field and
   $\{p_i:X_i\to Y_i:i\in I\}$  finitely many nonconstant morphisms of $k$-schemes of finite type.
  Then there are closed points  $x_i\in X_i$ with images $y_i=p_i(x_i)$ such that the residue fields $k(x_i), k(y_i)$ are all isomorphic.
  Moreover, the following hold.
  \begin{enumerate}
    \item If $k$ is infinite, then there is an  infinite set $D\subset \n$ such that, for every $d\in D$, there are
      infinitely many such points with $\deg(k(x_i)/k)=d$.
      \item 
  If $k$ is finite, then   there is an  infinite set $D\subset \n$ such that, for every $d\in D$, there are $\geq \frac12 |k|^d$ such points  with $\deg(k(x_i)/k)=d$.
\item If the $\{p_i:X_i\to Y_i:i\in I\}$ are chosen from a bounded family, and
  $\# I$ is bounded, then a fixed $D\subset \n$ works for all of them.
    \end{enumerate}
\end{lem}

Proof. We may assume that the $X_i, Y_i$ are affine and the $p_i$ are flat.
Choose dominant morphisms  $Y_i\to \a^1$ and let
$q:Z\to \a^1$ be the fiber product of all $X_i\to Y_i\to \a^1$.

Apply  \cite{MR1851662} to  $q:Z\to \a^1$.
We get a dense set of closed points    $z\in Z$ such that
$k(z)=k(q(z))$. Let $x_i\in X_i$ be the coordinate prjections of $z$. Then
$k(z)\supset k(x_i)\supset k(y_i)\supset k(q(z))$, so they are all equal.  

The claims (1--3) may be clearer from  \cite[7.2.1]{klos}. \qed

\medskip

A typical application is the following. 

\begin{cor}\label{same.res.field.cor.1}
  Let $k$ be  a perfect field, $X$ a $k$-scheme of finite type and
  $Z_1, Z_2\subset X$ positive dimensional, irreducible subsets.
  Let $$
  M\leftarrow C_M\stackrel{\pi}{\longrightarrow} X
  $$ be an irreducible family of curves joining $Z_1$ and $Z_2$. That is, 
$M\leftarrow C_M$ is flat, and there are 2  sections
  $\sigma_i: M\to C_M$ such that the $\pi\circ \sigma_i$ give dominant morphisms form $M$ to $Z_i$.
  Then  there is a dense set of closed points    $m\in M$ such that
  $k(m)=k\bigl(\pi\circ \sigma_1(m)\bigr)=k\bigl(\pi\circ \sigma_2(m)\bigr)$.\qed
  \end{cor}

%% \bibliography{../refs-main/refs}

\def\cprime{$'$} \def\cprime{$'$} \def\cprime{$'$} \def\cprime{$'$}
  \def\cprime{$'$} \def\dbar{\leavevmode\hbox to 0pt{\hskip.2ex
  \accent"16\hss}d} \def\cprime{$'$} \def\cprime{$'$}
  \def\polhk#1{\setbox0=\hbox{#1}{\ooalign{\hidewidth
  \lower1.5ex\hbox{`}\hidewidth\crcr\unhbox0}}} \def\cprime{$'$}
  \def\cprime{$'$} \def\cprime{$'$} \def\cprime{$'$}
  \def\polhk#1{\setbox0=\hbox{#1}{\ooalign{\hidewidth
  \lower1.5ex\hbox{`}\hidewidth\crcr\unhbox0}}} \def\cdprime{$''$}
  \def\cprime{$'$} \def\cprime{$'$} \def\cprime{$'$} \def\cprime{$'$}
\providecommand{\bysame}{\leavevmode\hbox to3em{\hrulefill}\thinspace}
\providecommand{\MR}{\relax\ifhmode\unskip\space\fi MR }
\providecommand{\MRhref}[2]{%
  \href{http://www.ams.org/mathscinet-getitem?mr=#1}{#2}
}
\providecommand{\href}[2]{#2}

  \bigskip

  Princeton University, Princeton NJ 08544-1000, \

  \email{kollar@math.princeton.edu}
\medskip

Beijing International Center for Mathematical Research,

Peking University, 100871, Beijing, China

  \email{zhiyutian@bicmr.pku.edu.cn}

\end{document}